\documentclass[11pt, a4paper, reqno]{amsart}

\usepackage{amssymb}
\usepackage{mathrsfs}
\usepackage{graphicx}

\usepackage[width=1.0\textwidth]{caption}

\usepackage[breaklinks]{hyperref}

\hypersetup{
unicode=true,
colorlinks=true,
linkcolor=blue,
citecolor=blue,
urlcolor=blue,
filecolor=blue,
bookmarksnumbered=true,
pdfstartview=FitH,
pdfhighlight=/N
}

\theoremstyle{plain}
\newtheorem{theorem}{Theorem}
\newtheorem{lemma}{Lemma}
\newtheorem{conj}{Conjecture}
\newtheorem{theoremSta}{Stahl Theorem}
\newtheorem{theoremBusTwoPoint}{Buslaev Two-Point Theorem}

\theoremstyle{definition}
\newtheorem{definition}{Definition}
\newtheorem{mcase}{Case}

\def\({\left(}
\def\){\right)}
\def\mcap{\operatorname{cap}}
\def\supp{\operatorname{supp}}
\def\const{\operatorname{const}}
\def\card{\operatorname{card}}
\def\bus{\operatorname{Bus}}
\def\Cheb{\operatorname{cheb}}

\def\NN{\mathbb N}
\def\RR{\mathbb R}
\def\CC{\mathbb C}

\def\ZZ{\mathbb Z}

\def\HH{\mathscr H}

\def\fH{\mathcal H}

\def\HH{\mathscr H}

\def\sL{\mathscr L}
\def\LL{\mathscr L}
\def\sN{\mathscr N}

\def\myA{\mathscr A}
\def\RS{\mathfrak R}
\def\myf{\mathfrak f}

\def\mypha{\vphantom{p}}

\let\leq\leqslant
\let\geq\geqslant

\let\eps\varepsilon

\let\myt\widetilde

\let\myo\overline

\begin{document}

\title{Zero Distribution of Hermite--Pad\'{e} Polynomials and Convergence
Properties of Hermite Approximants for Multivalued Analytic Functions}

\author{Nikolay R. Ikonomov}
\address{Institute of Mathematics and Informatics, Bulgarian Academy of Sciences}
\email{nikonomov@math.bas.bg}

\author{Ralitza K. Kovacheva}
\address{Institute of Mathematics and Informatics, Bulgarian Academy of Sciences}
\email{rkovach@math.bas.bg}

\author{Sergey P. Suetin}
\address{Steklov Mathematical Institute of Russian Academy of Sciences, Russia}
\email{suetin@mi.ras.ru}
\thanks{The third author was partially supported by the Russian Foundation
for Basic Research (RFBR, grants 13-01-12430-ofi-m2 and 15-01-07531-a),
and Russian Federation Presidential Program for support of the Leading
Scientific Schools (grant NSh-2900.2014.1).}

\date{March 10, 2016}

\begin{abstract}
In the paper, we propose two new conjectures about the convergence of Hermite
Approximants of multivalued analytic functions of Laguerre
class $\sL$. The conjectures are based in part on the numerical experiments,
made recently by the authors in \cite{KoIkSu15} and \cite{KoIkSu15b}.

Bibliography: \cite{Sze62} items.

Figures: \ref{Fig_la1sq(1500)320_full} items.
\end{abstract}

\maketitle


\markright{HERMITE--PAD\'{E} POLYNOMIALS}

\setcounter{tocdepth}{2}
\tableofcontents


\section{Introduction}\label{s1}

\subsection{Description of the problem}\label{s1s1}
The main goal of the current paper is to describe and illustrate the
main features of Hermite approximants of multivalued analytic
functions. The notion of Hermite approximants (HA) is very novel; it
was introduced in an implicit form by A.~Mart\'{\i}nez-Finkelshtein, E.~Rakhmanov, and
S.~Suetin in \cite{MaRaSu15}. We also propose two new conjectures (see
Conjecture \ref{conj1} and Conjecture \ref{conj2}) on convergence
properties of Hermite approximants of multivalued analytic functions
of Laguerre class $\sL$.
Given the germ of a function $f$ analytic at the point of infinity
$z=\infty$, the Hermite Approximants $\fH_{n,j}$, $j=0,1$, of order $n$ are completely determined by
the $3n+2$ initial Laurent coefficients of the given power series
of $f$. The rational functions $\fH_{n,j}, j=0,1$ are constructed on the basis
of type I Hermite--Pad\'e (HP) polynomials of the collection of three
functions\footnote{In what follows we suppose that the three functions
$1,f,f^2$ are rationally independent over the field of $\CC(z)$.}
$[1,f,f^2]$. All zeros and poles of HA are free. In this respect,
they are very similar to Pad\'e approximants (PA). From now on, we
assume that $f$ is a multivalued analytic function with a finite set of
singular points. In \cite{Sue15}, it was proven for a partial class of such multivalued
analytic functions that the HA is interpolating approximately $2n$ times (at free nodes)
some other branch\footnote{Partially, for the so-called
``differential-analytic functions''; about this notation
see \cite{ChCh84}, \cite{Hut06}.} of the given function $f$. Furthermore, there exist limit distributions
of the free zeros and poles of HA, as well as of the free nodes. The
associated limit measures solve some special equilibrium problems
for mixed Green logarithmic potentials with external fields. In some
particular cases, it was proven that such HA possesses an alternating
property which, as it turns out, is similar to the classical Chebysh\"ev's alternating
property. All these properties make HA very similar to the best
Chebysh\"ev rational approximants of analytic functions
(see \cite{GoRa87}, \cite{Rak15}). We note that for the construction of the
$n$-th HA $\fH_{n,j}$ merely the $3n+2$ initial Laurent coefficients suffice. In contrast
to this, in order to find the best Chebysh\"ev approximant one needs
the function $f$ to be given in an explicit form. Recall once again that
to construct the PA $[n/n]_f$ of order $n$ of the function $f$, given
by a power series, one should know the $2n+1$ initial Laurent
coefficients of the power series (see \cite{BaGr96}, \cite{ApBuMaSu11}). All
zeros and poles of PA
approximants are free, whereas the interpolation nodes are fixed at the point of infinity.
This approach is very novel and may be considered as a very
promising direction in the theory of the analytic continuation
(see \cite{Bib67}, \cite{Ara84}).

Let us now introduce the notion of HA of the analytic function $f$.
Given a germ $f$
\begin{equation}
\label{0.1}
f(z)=\sum_{k=0}^\infty\frac{c_k}{z^k}
\end{equation}
of a function $f$ analytic at the infinity point $z=\infty$, we assume
that the three functions
$1,f,f^2$ are rationally independent over the field of rational
functions $\CC(z)$. Let $n$ be fixed, $n\in\NN$.
Let now
$Q_{n,0}$, $Q_{n,1}$, $Q_{n,2}\in\CC^{*}_n[z]:= \CC_n[z]\setminus\{0\}$ be the
type I Hermite--Pad\'e (HP) polynomials of order ${n}$ for the collection
of the three functions $[1,f,f^2]$, that is
\begin{equation}
\label{0.2}
R_n(z):=(Q_{n,0}\cdot1+Q_{n,1}f+Q_{n,2}f^2)(z)=O\(\frac1{z^{2n+2}}\),
\quad z\to\infty.
\end{equation}
The polynomials $Q_{n,j}$ are not unique, but their ratios are
uniquely determined (see Lemma \ref{lem1} below). In what follows, we shall refer to the rational
functions $\fH_{n,0}:=Q_{n,0}/Q_{n,2}$ and $\fH_{n,1}:=Q_{n,1}/Q_{n,2}$ as
{\it Hermite approximants} of the given analytic function
$f\in\HH(\infty)$.
Given now a finite set $\Sigma\subset\CC$ (i.e. the set of finite cardinality,
$\card\Sigma<\infty$), we
denote by $\myA(\myo\CC\setminus\Sigma)$ the class of all functions
$f\in\HH(\infty)$ which admit an analytical continuation from the infinity point along each path
 avoiding the given set $\Sigma$. Let
$\myA^\circ(\myo\CC\setminus\Sigma):=\myA(\myo\CC\setminus\Sigma)\setminus
\HH(\myo\CC\setminus\Sigma)$. Up to the end of the current paper, we will restrict our attention, while
discussing problems concerned with HA, only to the
{\it Laguerre class} $\sL$ of multivalued analytic functions. In other words, to the
class of multivalued analytic functions given by the explicit
representation
\begin{equation}
\label{l1}
f(z)=\prod_{j=1}^p(z-a_j)^{\alpha_j},
\quad\alpha_j\in\CC\setminus\ZZ,\quad\sum_{j=1}^p\alpha_j=0,
\end{equation}
where $a_j\in\CC$, $j=1,\dots,p$, and $a_j\neq a_k,\,
j\neq k$. Thus if $f\in\LL$, then
$f\in\myA^\circ(\myo\CC\setminus\Sigma)$, where $\Sigma=\{a_1,\dots,a_p\}$.
Let us fix the germ of $f$ at the infinity point
by the condition $f(\infty)=1$.

We mainly restrict our attention to the partial subclass $\LL$ of $\sL_{\RR}$
 of functions
 given by the representation
\begin{equation}
\label{l2}
f(z)=\prod_{j=1}^q\(\frac{z-e_{2j-1}}{z-e_{2j}}\)^{\alpha_j},
\quad\alpha_j\in\RR\setminus\ZZ,
\end{equation}
where $e_j\in\RR$, $j=1,\dots,2q$, and $e_1<\dots<e_{2q}$.

\subsection{New conjectures}\label{s1s2}
The main purpose of the current paper is to explain how to use the
{\it Hermite Approximants} (HA) in the constructive approximation theory, as well as
to impose two new conjectures on the convergence of HA for the functions of
Laguerre class $\sL$. To be more precise, we are interested in studying
type I Hermite--Pad\'e polynomials and the corresponding rational HA
with free zeros and poles, as well as interpolation nodes.

The main objectives of the current paper are the next conjectures.

\begin{conj}\label{conj1}
Let $f\in\LL$ and the functions $1,f,f^2$ be rationally independent
over the field $\CC(z)$. Then for $z\in\myo\CC\setminus F$
\begin{equation}
\label{co1}
\frac{Q_{n,0}}{Q_{n,2}}(z)\overset{\mcap}\longrightarrow f^2(z),\quad n\to\infty,
\end{equation}
where the compact set $F$ consists of a finite number of closed analytic arcs,
$F=\bigcup_{j=1}^m F_j$.
\end{conj}

\begin{conj}\label{conj2}
Let $f\in\LL$ and all the exponents $\alpha_j\neq\pm1/2$ (see \eqref{l1}).
Then for $z\in\myo\CC\setminus F$
\begin{equation}
\label{co2}
\frac{Q_{n,1}}{Q_{n,2}}(z)\overset{\mcap}\longrightarrow\const \cdot f(z),\quad n\to\infty,
\end{equation}
where $\const\neq0$ and the compact set $F$ is just the same as in
Conjecture \ref{conj1}.
\end{conj}

Conjectures \ref{conj1} and \ref{conj2} might be considered as a step
towards the construction of a general convergence theory of Hermite
Approximants. No doubt that the new theory should be much more complicated than
Stahl's Theory about classical PA and Buslaev's Theory about multipoint PA.
For other conjectures on the limit zero distribution (LZD)
of HP polynomials, the reader is referred to  \cite{Nut84}, \cite{Sta88} and \cite{Apt08}.

\begin{figure}[!ht]
\centerline{
\includegraphics[width=15cm,height=15cm]{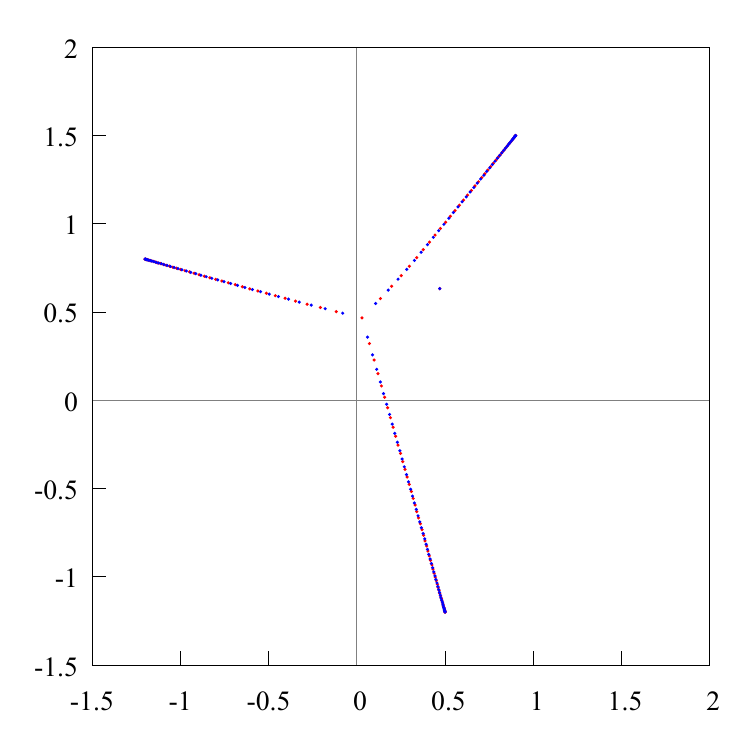}}
\vskip-6mm
\caption{Zeros (blue points) and poles (red points) of PA $[130/130]_f$ of
the function
$f(z)=\bigl(z-(-1.2+0.8i)\bigr)^{1/3}\bigl(z-(0.9+1.5i)\bigr)^{1/3}
\bigl(z-(0.5-1.2i)\bigr)^{-2/3}$. Since the genus of the corresponding Stahl's
two-sheeted Riemann surface equals $1$ (i.e. it is an elliptic Riemann
surface), there might be at most a single ``spurious'' zero-pole pair, i.e.
a single Froissart doublet. It is really present on the picture.
}
\label{Fig_pade10_2500_130}
\end{figure}

\begin{figure}[!ht]
\centerline{
\includegraphics[width=15cm,height=15cm]{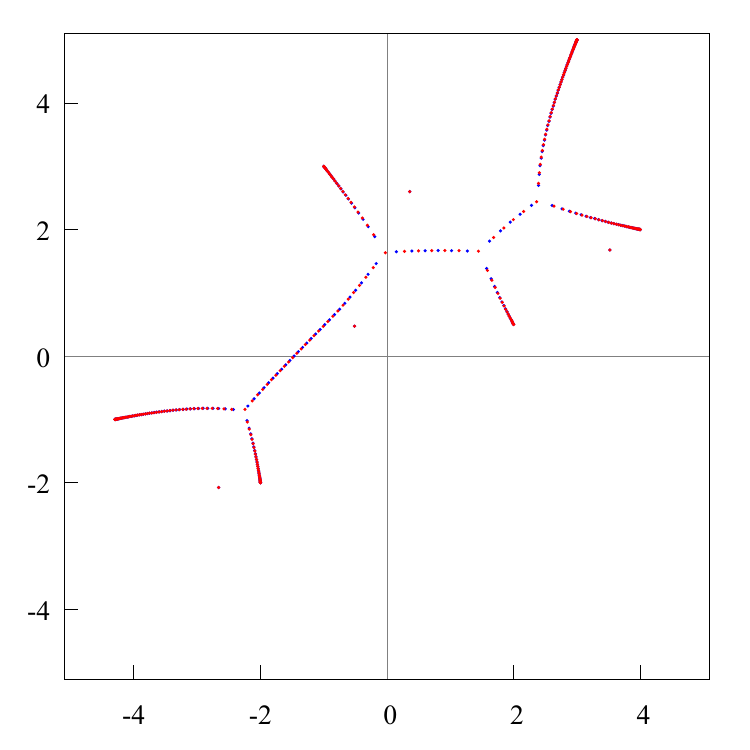}}
\vskip-6mm
\caption{Zeros (blue points) and poles (red points) of PA $[267/267]_f$ of
the function
$f(z)=\{(z+(4.3+1.0i))
(z-(2.0+0.5i))(z+(2.0+2.0i))(z+(1.0-3.0i))
(z-(4.0+2.0i))(z-(3.0+5.0i))\}^{-1/6}$.
These zeros and poles are distributed in a plane, under fixed $n=267$,
accordingly to the electrostatic model by Rakhmanov \cite{Rak12}.
There are $4$ Chebotar\"ev points on the picture. Thus the genus
of the Stahl's hyperelliptic Riemann surface is $4$.
By this reason for each $n$ there might be no more than $4$ Froissart doublets. Here
are observed $4$ Froissart doublets (cf. \cite[Fig. 2]{Sue15}).
}
\label{Fig_pade103_5000_267_full}
\end{figure}

\begin{figure}[!ht]
\centerline{
\includegraphics[width=15cm,height=15cm]{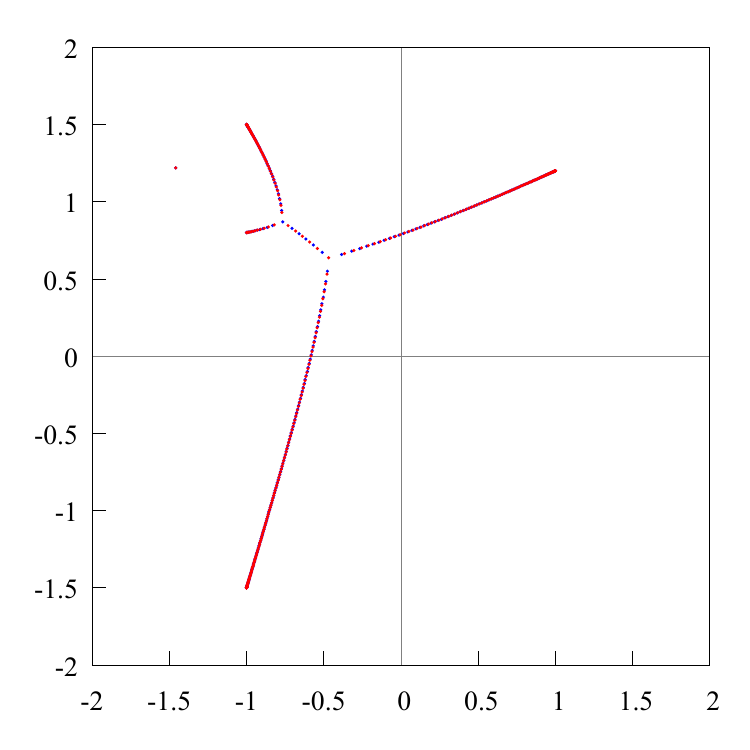}}
\vskip-6mm
\caption{Zeros (blue points) and poles (red points) of PA $[300/300]_f$ of
the quadratic function
$f(z)=\(\dfrac{z-(-1.0+0.8i)}{z-(1.0+1.2i)}\)^{1/2}
+\(\dfrac{z-(-1.0+1.5i)}{z-(-1.0-1.5i)}\)^{1/2}$.
There are $2$ Chebotar\"ev points on the picture. Thus the genus
of the Stahl's hyperelliptic Riemann surface is $2$.
Here we observe a single Froissart doublet located in the second quadrant.
In full compliance with the Rakhmanov's model \cite{Rak12}, the Froissart
doublet attracts to itself the Stahl's $S$-compact set $S_{300}$;
cf. Fig. \ref{Fig_PA_t_log_1_300(2000)_full}.
}
\label{Fig_Pade_6(3000)300_full}
\end{figure}

\begin{figure}[!ht]
\centerline{
\includegraphics[width=15cm,height=15cm]{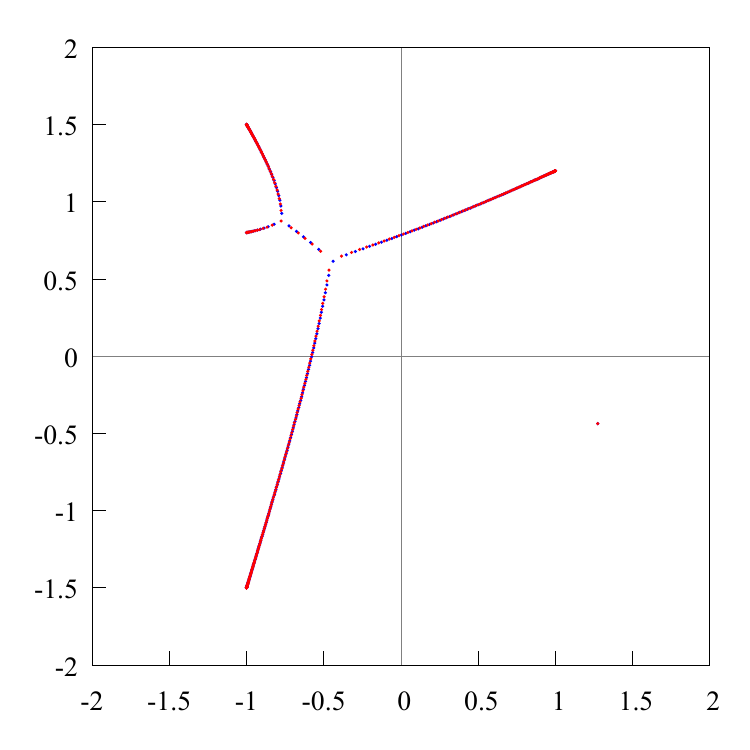}}
\vskip-6mm
\caption{Zeros (blue points) and poles (red points) of PA $[300/300]_f$ of
the logarithmic function
$f(z)=\log\(\dfrac{z-(-1.0+0.8i)}{z-(1.0+1.2i)}\)
+\log\(\dfrac{z-(-1.0+1.5i)}{z-(-1.0-1.5i)}\)$.
There are $2$ Chebotar\"ev points on the picture. Thus the genus
of the Stahl's hyperelliptic Riemann surface is $2$.
Here we observe a single Froissart doublet located in the fourth quadrant.
In full compliance with the Rakhmanov's model \cite{Rak12}, the Froissart
doublet attracts to itself the Stahl's $S$-compact set $S_{300}$;
cf. Fig. \ref{Fig_Pade_6(3000)300_full}.
}
\label{Fig_PA_t_log_1_300(2000)_full}
\end{figure}

\begin{figure}[!ht]
\centerline{
\includegraphics[width=15cm,height=15cm]{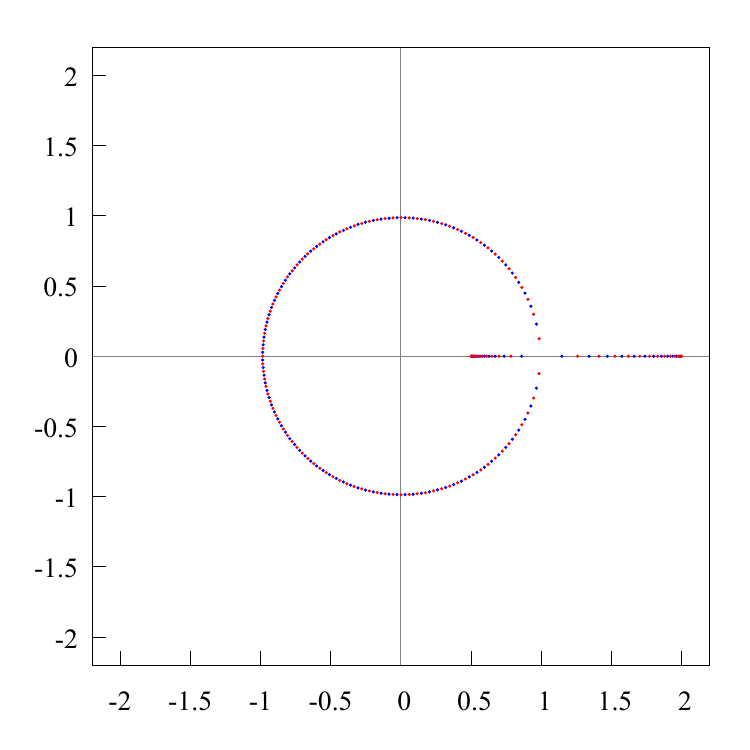}}
\vskip-6mm
\caption{Numerical zeros (blue points) and poles (red points) distribution of
two-point PA $[120/120]_{\myf}$ of the set of functions $\myf=\{f_0,f_\infty\}$,
where $f_0=((1-2z)(2-z))^{-1/2}$, $f_0\in\HH(0)$,
$f_\infty=((2z-1)(z-2))^{-1/2}+1$, $f_\infty\in\HH(\infty)$. The germs $f_0$ and $f_\infty$ result in
two different multivalued analytic functions. Thus, this is a generic case
and by Buslaev's Theorem the associated weighted $S$-curve divides
the Riemann sphere into two domains.}
\label{Fig_bus210b_4000_120_full}
\end{figure}

\begin{figure}[!ht]
\centerline{
\includegraphics[width=15cm,height=15cm]{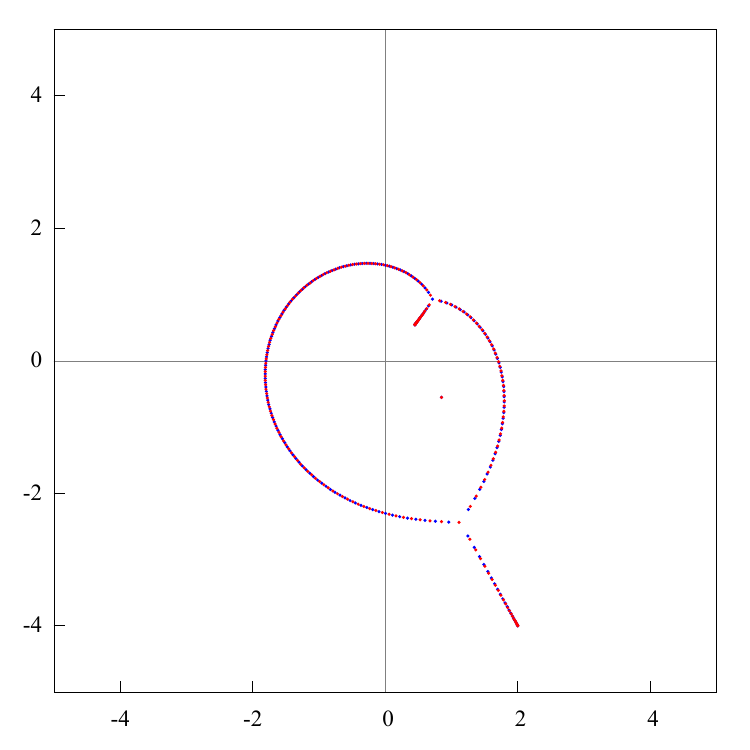}}
\vskip-6mm
\caption{Numerical zeros (blue points) and poles (red points) distribution of
two-point PA $[195/195]_f$ to the function $f(z)=\sqrt[4]{(z-a_1)/(z-a_2)}$,
where $a_1=0.9-1.1i$ and $a_2=0.1+0.2i$.
Here are selected two ``quite different branches''
of the function $f$, namely,
$f_0=\sqrt[4]{(z-a_1)/(z-a_2)}$ and $f_\infty=-\sqrt[4]{(z-a_1)/(z-a_2)}$.
All, but one pair, zeros (blue points) and poles (red points) approximate
numerically Buslaev's compact set. But there is a single Froissart doublet
located in the domain $D_0(f)\ni0$; cf. \cite[Fig. 3]{Sue15}.
}
\label{Fig_bus205c(2000)195_full}
\end{figure}

\begin{figure}[!ht]
\centerline{
\includegraphics[width=15cm,height=15cm]{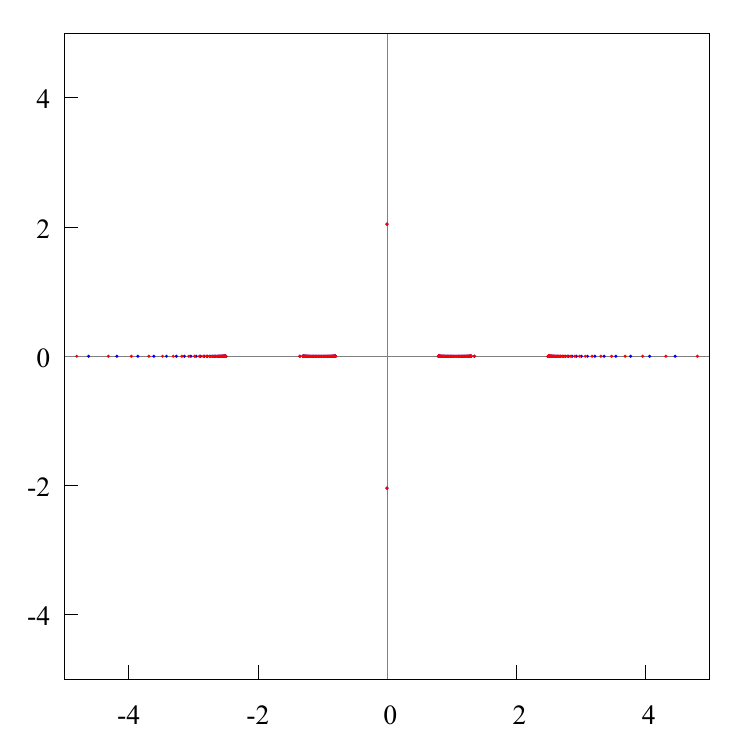}}
\vskip-6mm
\caption{Zeros of HP polynomials $Q_{200,0}$ (blue points) and
$Q_{200,1}$ (red points) for the triple of functions $[1,f,f^2]$, where
$f(z)=
\(\dfrac{z+2.5}{z+1.3}\)^{1/3}\(\dfrac{z+0.8}{z-0.8}\)^{1/3}
\(\dfrac{z-1.3}{z-2.5}\)^{1/3}$.
All but two pairs of zeros are distributed in accordance with
Theorem \ref{the1} on the real line on the complement of the three closed
segments $[-2.5,-1.3]$, $[-0.8,0.8]$, and $[1.3,2.5]$. There are two pairs
of complex conjugate Froissart doublets; cf. Fig. \ref{Fig_La1(2000)200_full}.
}
\label{Fig_La1(2000)200_blu_red}
\end{figure}

\begin{figure}[!ht]
\centerline{
\includegraphics[width=15cm,height=15cm]{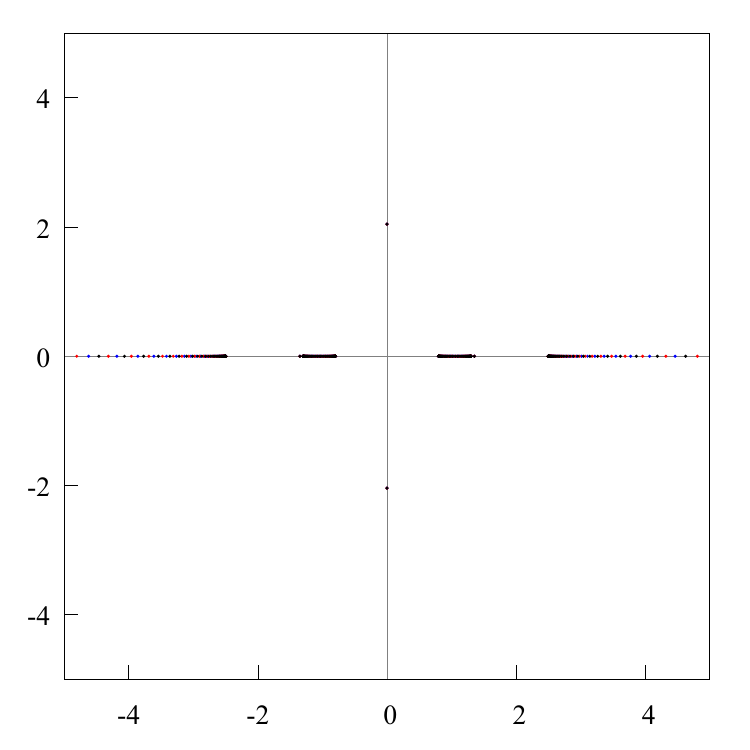}}
\vskip-6mm
\caption{Zeros of HP polynomials $Q_{200,0}$ (blue points), $Q_{200,1}$ (red
points), and $Q_{200,2}$ (black points) for the triple of functions
$[1,f,f^2]$, where
$f(z)=
\(\dfrac{z+2.5}{z+1.3}\)^{1/3}\(\dfrac{z+0.8}{z-0.8}\)^{1/3}
\(\dfrac{z-1.3}{z-2.5}\)^{1/3}$.
All but two pairs of zeros are distributed in accordance with
Theorem \ref{the1} on the real line on the complement of the three closed
segments $[-2.5,-1.3]$, $[-0.8,0.8]$, and $[1.3,2.5]$. There are two pairs
of complex conjugate Froissart triplets; cf. Fig. \ref{Fig_La1(2000)200_blu_red}.
}
\label{Fig_La1(2000)200_full}
\end{figure}

Conjectures \ref{conj1} and \ref{conj2} are based on the rigorous results
from \cite[Theorem 1.8]{MaRaSu15} and \cite{Sue15}, as well as on the numerical
experiments produced by the authors
in \cite{KoIkSu15}, \cite{KoIkSu15b}; for more details see \S \ref{s3s4}
and Fig. \ref{Fig_pade10_2500_130}--\ref{Fig_la1sq(1500)320_full}.

\section{Pad\'e approximants}\label{s2}

\subsection{Pad\'e approximants and \texorpdfstring{$J$}{J}-fractions}\label{s2s1}
Since Hermite approximants are a generalization of classical Pad\'e
approximants, we start from the basic definition of Pad\'e polynomials
$P_{n,0},P_{n,1}$ and of
Pad\'e approximants $[n/n]_f:=-P_{n,0}/P_{n,1}$.

We recall the definition of PA of an analytic function, given by the power
(in fact, by Laurent) series \eqref{0.1} at the infinity point $z=\infty$. For the sake of
 convenience, we introduce the Pad\'e polynomials
$P_{n,0},P_{n,1}\in\CC^{*}_n[z]$ in the following way. There exist
 polynomials $P_{n,0}$ and $P_{n,1}$ of degree $\leq{n}$ such that (cf. \eqref{0.2})
\begin{equation}
\label{1.4p}
(P_{n,0}\cdot1+P_{n,1}\cdot f)(z)=O\(\frac1{z^{n+1}}\),
\quad z\to\infty.
\end{equation}
The rational function $-P_{n,0}/P_{n,1}$ is uniquely determined and is called
the diagonal Pad\'e approximant $[n/n]_f=-P_{n,0}/P_{n,1}$ of the function $f$
(at the infinity point).

In the ``generic case'' the relation \eqref{1.4p} is equivalent to the relation
\begin{equation}
\label{1.5p}
\bigl(f-[n/n]\bigr)(z)=O\(\frac1{z^{2n+1}}\),\quad z\to\infty.
\end{equation}
Thus, from \eqref{1.5p} follows that the $n$-th PA $[n/n]_f$ is the {\it best
local rational} approximant of order $\leq{n}$ of the given power
series \eqref{0.1}. Notice that $[n/n]_f$ is a rational function with free poles
and free zeros. Furthermore, it interpolates the given power series $f$ at the
fixed point $z=\infty$ up to the order $2n+1$. Hence,
\begin{equation}
\label{1.5.2}
[n/n]_f(z)=c_0+\frac{c_1}z+\dots+\frac{c_{2n}}{z^{2n}}+O\(\frac1{z^{2n+1}}\),
\quad z\to\infty.
\end{equation}
Recall that by definition of the partial sums $S_{2n}(z)$ of the power
series \eqref{0.1}
\begin{equation}
\label{1.5.3}
S_{2n}(z)=c_0+\frac{c_1}z+\dots+\frac{c_{2n}}{z^{2n}},
\end{equation}
that is,
$$
[n/n]_f(z)-S_{2n}(z)=O\(\frac1{z^{2n+1}}\),
\quad z\to\infty.
$$
Relations  \eqref{1.5.2} and \eqref{1.5.3} together lead to a very natural question, namely: do the PA
$[n/n]_f(z)$ have some real advantages over the partial sums $S_{2n}(z)$?

The answer is ``yes'' and comes from the classical $J$-fractions theory.
This is the well-known classical way to evaluate an analytic function going out
from its germ, a way
which goes back to Gauss and Jacoby (see \cite{BaGr96}).
However, they used it to evaluate only special (in particular, hypergeometric) functions.

Recall that $f$ is a multivalued analytic function with a
finite set $\Sigma$ of branch points, $\card\Sigma<\infty$. To be more
precise, we suppose that $f$ is analytic in the domain
$\myo\CC\setminus\Sigma$, but not holomorphic in
$\myo\CC\setminus\Sigma$. We adopt the notation
$f\in\myA^\circ(\myo\CC\setminus\Sigma):
=\myA(\myo\CC\setminus\Sigma)\setminus\HH(\myo\CC\setminus\Sigma)$.

Let $Q\in\CC[z]$ be an arbitrary complex polynomial.
We denote the zero-counting measure of the polynomial $Q$ by $\chi(Q)$,
that is,
\begin{equation}
\label{zero-co}
\chi(Q):=\sum_{\zeta:Q(\zeta)=0}\delta_\zeta,
\end{equation}
where the zeros of $Q$ are counted with regards to their multiplicities; as
usual, $\delta_\zeta$ denotes the Dirac measure, concentrated at the point
$\zeta\in\myo\CC$.

Let $f\in\sL$ be in a ``generic case''.

Then we can use the functional analog of
Euclid's algorithm to obtain the formal expansion (see \cite{Che55}, \cite{BaGr96})
\begin{align}
f(z)&=1+\frac{A_1}{z-B_1-f_1(z)}
=1+\frac{A_1}{z-B_1-\dfrac{A_2}{z-B_2-f_2(z)}}\notag\\
&\simeq1+{\frac{{A_1}\vphantom b\hfill}
{
z-B_1-\dfrac{A_2\vphantom b\hfill}
{z-B_2-\dfrac{A_3\qquad\vphantom b\hfill}{\smallmatrix\\\ddots\endsmallmatrix}
}}}\simeq J(z),
\notag
\end{align}
where all $A_n$ do not vanish, $A_n\neq 0$, $n=2,3,\dotsc$. As usual, the notation
``$\simeq$'' means only a formal equality with no
convergence statements.
Thus let us consider the $n$-th truncate $J_n(z)$ of the continued fraction
$J(z)$, i.e.
$$
J_n(z):=1+\frac{A_1\vphantom b\hfill}
{z-B_1-\dfrac{A_2\vphantom b\hfill}
{z-B_2-\dfrac{A_3\qquad\vphantom b\hfill}
{\dfrac{\qquad\ddots\qquad}
{z-B_{n-1}-\dfrac{A_n}{z-B_n}}.
}}}
\notag
$$
We recall that $J_n$ is a rational function of order $n$, $J_n\in\CC_n(z)$. Set
$J_\infty(z):=\lim\limits_{n\to\infty}J_n(z)$. The problem of
convergence of the $J$-fraction to $f$ may be stated as the problem of
equality
\begin{equation}
\label{1.3}
f(z)\overset{\text{?}}=J_\infty(z),
\end{equation}
in other words, it is the problem of evaluation of $f(z)$ via $J_\infty(z)$.
Since $f$ is a multivalued function and all the $J_n$ are single valued
functions, two main questions arise in connection with
Problem \eqref{1.3}: {\it in what sense this equality might be understood
and in what domain does it hold true}?

\subsection{Problem of equality \texorpdfstring{$f(z)=J_\infty(z)$: case $p=2$}{f(z)=Jinfty(z): case p=2}}\label{s2s2}
Let in \eqref{l1} $p=2$, i.e.
\begin{equation}
\label{1.32}
f(z)=\(\frac{z+1}{z-1}\)^\alpha,
\quad \alpha\in\CC\setminus\ZZ.
\end{equation}
Let assume that $\alpha\in(-1/2,1/2)$, $\alpha\neq0$. Then
for $Q_n$, where $J_n=P_n/Q_n$, we easily obtain (see \cite{Che55}, \cite{Sze62})
\begin{equation}
\label{ja1}
\int_{-1}^1 Q_n(x)\,x^k\(\frac{1+x}{1-x}\)^\alpha\,dx=0,\quad k=0,\dots,n-1.
\end{equation}
Thus $Q_n(x)=P^{(-\alpha,\alpha)}_n(x)$ is the Jacobi polynomial of
degree $n$ with the parameters $(-\alpha,\alpha)$, $\alpha\in(-1/2,1/2)$,
and orthogonal on $E=[-1,1]$. It follows  immediately from  \eqref{ja1}
 that all zeros of $Q_n$ belong to the segment
$[-1,1]$, furthermore, for $\chi(Q_n)$ the
relation
\begin{equation}
\frac1n\chi(Q_n)
\overset{*}\longrightarrow\frac{dx}{\pi\sqrt{1-x^2}},
\quad n\to\infty.
\end{equation}
holds (see \eqref{zero-co}).
It is well-known that $P^{(-\alpha,\alpha)}_n(z)$ solves
the following linear differential equation of 2-nd degree
\begin{equation}
(z^2-1)w''+2(z-\alpha)w'-n(n+1)w=0.
\label{ja2}
\end{equation}
By applying the classical
asymptotic Liouville--Steklov method \cite[\S~8.63]{Sze62} to the equation \eqref{ja2},
we obtain a formula for the strong asymptotics of the Jacobi polynomials
\begin{equation}
P_n^{(-\alpha,\alpha)}(z)=
\(\frac{z-1}{z+1}\)^{\alpha/2}\frac{(z+(z^2-1)^{1/2})^{n+1/2}}{(z^2-1)^{1/4}}
\(1+O\(\frac1n\)\),\quad z\notin\Delta.
\label{jac1}
\end{equation}
Since for the numerator $P_n$ of $J_n=P_n/Q_n$ we have
$P_n=P^{(\alpha,-\alpha)}_n$, i.e. $P_n$ is Jacobi polynomial of order
$n$ with parameters $(\alpha,-\alpha)$, a direct analog
of the strong asymptotics formula \eqref{jac1} is also valid for $P_n$.
Thus after combining them, these two formulae provide a strong
asymptotics formula for the rational function $J_n$
$$
J_n(z)=\(\frac{z+1}{z-1}\)^\alpha\(1+O\(\frac1n\)\),\quad z\notin\Delta.
$$
Therefore,
\begin{equation}
\label{jac_eq}
f(z)=J_\infty(z)\quad\text{for}\quad z\in\myo\CC\setminus\Delta.
\end{equation}

\subsection{Problem of equality \texorpdfstring{$f(z)=J_\infty(z)$: case $p=3$}{f(z)=Jinfty(z): case p=3}}\label{s2s3}
In 1885 Laguerre \cite{Lag85} made an
attempt to solve Problem \eqref{1.3} for the partial case when $p=3$ and $f\in\LL$
(cf. \eqref{l1}), i.e.,
\begin{equation}
\label{1.4}
f(z)=\prod_{j=1}^3(z-a_j)^{\alpha_j},
\quad \alpha_j\in\CC\setminus\ZZ,\quad\sum_{j=1}^3\alpha_j=0,
\end{equation}
where the points $a_1,a_2,a_3$ are in a ``general position''; in particular, they
are pairwise distinct and don't belong to a straight line.

Laguerre derived in 1885 the property of nonhermittian
orthogonality for the denominators $Q_n$ of the rational function
$J_n=P_n/Q_n$, i.e.
\begin{equation}
\oint_\Gamma Q_n(\zeta)\zeta^k f(\zeta)\,d\zeta=0,\quad k=0,\dots,n-1,
\label{1.5}
\end{equation}
where $\Gamma$ is an arbitrary closed contour that separates the three
points $a_1$, $a_2$, $a_3$ from the infinity point. He also proved (see
also \cite{Per57}, \cite{Nut86}, \cite{MaRaSu15}) that the
polynomial $P_n$ and the function $Q_nf$ solve the following linear
differential equation of second order
\begin{equation}
\label{lag1}
A_3(z)\Pi_{n,1}(z)w''+\Pi_{n,3}(z)w'+\Pi_{n,2}(z)w=0,
\end{equation}
where $A_3(z)=\prod_{j=1}^3(z-a_j)$ and $\Pi_{n,k}\in\CC_k[z]$,
$k=1,2,3$, are some polynomials of degree $k$. To be more precise,
\begin{align*}
\Pi_{n,1}(z)&=z-z_n , \\
\Pi_{n,2}(z)&=-n(n+1)(z-b_n)(z-v_n) , \\
\Pi_{n,3}(z)&=(z-z_n)B_2(z)-A_3(z) ,\quad B_2(z)=A_3'(z)f'(z)/f(z) .
\end{align*}
Thus the polynomial coefficients in equation \eqref{lag1} are of fixed
degrees, but depend on $n$. These polynomial
coefficients contain three so-called accessory parameters $z_n,b_n,v_n$,
the behavior of which as $n\to\infty$ is presently unknown. That is why
Laguerre couldn't solve neither the problem about the asymptotic behavior of the
polynomials $P_n$ and $Q_n$, nor the Problem about the equality
$f(z)=J_\infty(z)$ as well.

For case $\Sigma=\{a_1,a_2,a_3\}$, Problem \eqref{1.3},
which is about the strong convergence of $J$-fraction,
was solved by J. Nuttall in 1986 only in terms
of PA, and on the basis of the seminal Stahl's Theorem \cite{Sta97b} about the convergence
in capacity of PA of an arbitrary multivalued analytic function with a finite set of branch
points (for the strong asymptotics and strong convergence properties, see also \cite{Nut84},
\cite{Nut90}, \cite{Sue00}, \cite{BaYa09}, \cite{KoIkSu15c}, \cite{ApYa11}).

In the ``generic case'' $[n/n]_f(z)=J_n(z)$. Hence, the
Problem about the equality $f(z)=J_\infty(z)$ is in fact the problem about the
convergence of the sequence of PA $\{[n/n]_f(z),n=0,1,\dots\}$ of the given
analytic function $f$.

Nuttall proved (see  \cite{Nut86}) that for the function $f\in\sL$, given
by \eqref{1.4}, the equality $f(z)=J_\infty(z)$ holds true inside the
domain $D:=\myo\CC\setminus{S}$, where $S$ is Stahl's compact set,
up to a unique arbitrary zero-pole pair (in other words, a
spurious zero-pole pair, or a Froissart doublet;
see \cite{Fro69}, \cite{Sue00}, \cite{ApYa11}). To be more precise,
there is a sequence $z_n\in\myo\CC$ such that for each compact set
$K\subset{D}$ and for every positive $\eps>0$
\begin{equation}
\label{nut86}
\sup_{z\in K\setminus\{z:|z-z_n|<\eps\}}\bigl|f(z)-J_n(z)\bigr|\to0,\quad
n\to\infty
\end{equation}
 (cf. \eqref{jac_eq}).
Notice that the convergence relation \eqref{nut86} does not result from Stahl's Theorem,
since it is dealing with the LZD (the Limit Zero Distribution) of Pad\'e
polynomials and with the convergence of PA in capacity; for the strong
convergence see also \cite{BaYa09}, \cite{ApYa11}, \cite{MaRaSu12}, \cite{KoSu14}.


\subsection{Classical Pad\'e approximants: Stahl's Theory}\label{s2s4}
Let $f\in\HH(\infty)$ be a multivalued analytic function in the class
$\myA^\circ(\myo\CC\setminus\Sigma)$, $\card\Sigma<\infty$,

Given a positive Borel measure $\mu$ with a compact support
$\supp\mu\subset\myo\CC$, $\supp\mu\neq\myo\CC$, let
$V^\mu(z)$ be the logarithmic potential (see \cite{Lan66}, \cite{SaTo97})
associated with $\mu$, that is:
$$
V^\mu(z):=
\int_{\supp\mu}\log\frac1{|z-\zeta|}\,d\mu(\zeta).
$$
We set $V_{*}^\mu(z)$ for the spherically normalized
logarithmic potential of measure $\mu$, i.e.
$$
V_{*}^\mu(z):=
\int_{|\zeta|\leq1}\log\frac1{|z-\zeta|}\,d\mu(\zeta)
+\int_{|\zeta|>1}\log\frac1{|1-z/\zeta|}\,d\mu(\zeta).
$$

Let $f$ be the germ of a
multivalued analytic function $f$ with a finite set of branch
points. Then the seminal Stahl's Theorem
gives a complete answer to the problem about the limit zero-pole
distribution of the classical PA~of~$f$.
The keystone of Stahl's Theory is the existence of a unique
``maximal domain'' of holomorphy of $f$, i.e. of a domain
$D=D(f)\ni\infty$ such that the given germ $f$ can be continued as a
holomorphic (i.e. analytic and single-valued) function from a
neighborhood of the infinity point $z=\infty$ into~$D$ (i.e. the function $f$
is continued analytically along each path belonging to~$D$). ``Maximal''
means that
$\partial{D}$ is of ``minimal capacity'' among all compact sets
$\partial{G}$ such that $G$ is a domain, $G\ni\infty$ and $f\in\HH(G)$.
To be more precise, we have
$$
\mcap{\partial D}
=\min\{\mcap{\partial G}:\,\text{domain}\,G\ni\infty,f\in\HH(G)\}.
$$
The ``maximal'' domain $D$ is unique up to an arbitrary compact set of zero
capacity.

The compact set $S=S(f):=\partial{D}$ is called ``Stahl's
compact set'' or ``Stahl's $S$-compact set'' and $D$ is called ``Stahl's
domain'', respectively. The crucial properties of $S$ for the theory of Stahl to be
true are the following: the complement $D=\myo\CC\setminus{S}$ is a
domain, $S$ consists of a finite number of analytic arcs (in fact, the
union of the closures of the critical trajectories of a quadratic
differential), and finally, $S$ possesses the following property of
``symmetry'' (compact sets of such type are usually called
 ``$S$-compact sets'' or ``$S$-curves'', see \cite{Rak12}, \cite{KuSi15}):
\begin{equation}
\frac{\partial g^{}_S(z,\infty)}{\partial n^{+}}
=\frac{\partial g^{}_S(z,\infty)}{\partial n^{-}},\quad z\in S^\circ;
\label{1.6}
\end{equation}
where $g_S(z,\infty)$ is the Green's function of the domain $D$
with a logarithmic singularity at the point $z=\infty$, $S^\circ$ is the
union of all open arcs of $S$ (whose closures constitute $S$, i.e.
$S\setminus S^\circ$ is a finite set), and $\partial n^{+},\,\partial
n^{-}$ mean the inner (with respect to $D$) normal derivatives of
$g_S(z,\infty)$ at a point $z\in S^\circ$ from the opposite sides of $S^\circ$.
Let $\lambda=\lambda_S$ be the unique equilibrium probability measure for $S$,
i.e. $V^\lambda(z)\equiv\const=\gamma_S$ for $z\in S$; $\gamma_S$ is
the Robin constant for $S$. Then, by the identity
$g_S(z,\infty)\equiv\gamma_S-V^\lambda(z)$, the property of
symmetry \eqref{1.6} is equivalent to the property
\begin{equation}
\frac{\partial V^\lambda}{\partial n^{+}}(z)
=\frac{\partial V^\lambda}{\partial n^{-}}(z),\quad z\in S^\circ.
\label{1.6.1}
\end{equation}
If
$$
f(z)=\prod_{j=1}^3(z-a_j)^{\alpha_j},
$$
then we have that the compact set $S$ consists of the critical trajectories of the
quadratic differential
\begin{equation}
-\frac{z-v}{A_3(z)}\,dz^2>0,\quad A_3(z):=\prod_{j=1}^3(z-a_j).
\label{q2}
\end{equation}
These trajectories emanate from the points $a_j$ and culminate at the so-called
Chebotar\"ev's point $z=v$ (see \cite{Kuz80}). All points $a_1,a_2,a_3$ are the simple
poles of the quadratic differential \eqref{q2} and the Chebotar\"ev point
is the simple zero of that differential. In general, Chebotar\"ev's point couldn't be found
via elementary functions of the points $a_1,a_2,a_3$.
It is uniquely determined from
the condition that both periods of the Abelian integral
\begin{equation}
\int^z\sqrt{\frac{\zeta-v}{A_3(\zeta)}}\,d\zeta
\label{a3}
\end{equation}
are purely imaginary.
Because of this, the function
\begin{equation}
\Re
\int_{a_1}^z\sqrt{\frac{\zeta-v}{A_3(\zeta)}}\,d\zeta
\label{3.14}
\end{equation}
is a single-valued harmonic function on the two-sheeted elliptic
Riemann surface $\RS_2$, given by the equation $w^2=(z-v)A_3(z)$. The
Chebotar\"ev--Stahl compact set is given by the equality
\begin{equation}
S=\biggl\{z\in\CC:
\Re
\int_{a_1}^z\sqrt{\frac{\zeta-v}{A_3(\zeta)}}\,d\zeta=0
\biggr\},
\label{3.15}
\end{equation}
and the so-called $g$-function
\begin{equation}
g(z):=\Re
\int_{a_1}^z\sqrt{\frac{\zeta-v}{A_3(\zeta)}}\,d\zeta
\label{3.16}
\end{equation}
equals identically to the Green's function $g_S(z,\infty)$ of the
domain $D$. From the above results it follows immediately that
for the equilibrium measure $\lambda$ (see~\eqref{q2}) the following
representation holds:
$$
d\lambda(z)=\frac1{\pi i}\sqrt{\frac{z-v}{A_3(z)}}\,dz>0,\quad z\in S.
$$

One of the main results of Stahl's Theory (see \cite{Sta85a}--\cite{Sta86b}, and also \cite{Sta97b}) is

\begin{theoremSta}[H. Stahl, 1985--1986]
Let the function $f\in\HH(\infty)$, $f\in\myA^\circ(\myo\CC\setminus\Sigma$),
$\card\Sigma<\infty$, let $D=D(f)$ be Stahl's ``maximal'' domain of $f$,
$S=\partial{D}$ be Stahl's compact set, and
$[n/n]_{f}=-P_{n,0}/P_{n,1}$ be the $n$-th diagonal Pad\'e approximant of
the function $f$. Then the following statements are valid:

1) There exists a LZD of Pad\'e polynomials $P_{n,j}$,
$j=0,1$, namely,
\begin{equation}
\frac1n\chi(P_{n,j})
\overset{*}\longrightarrow\lambda,\quad\text{as}\quad
n\to\infty,\quad j=0,1,
\label{2.1}
\end{equation}
where $\lambda=\lambda_S$ is the unique probability equilibrium
measure for the compact set $S$, i.e. $V^\lambda(z)\equiv\gamma_S$,
$z\in S$, $\gamma_S$ -- the Robin constant for $S$;

2) the $n$-th diagonal Pad\'e approximants converge in capacity to the
function $f$ inside the domain $D$,
\begin{equation}
[n/n]_f(z)\overset{\mcap}\longrightarrow f(z),\quad n\to\infty,\quad z\in
D;
\label{2.2}
\end{equation}

3) the rate of the convergence in \eqref{2.2} is completely characterized by the equality
\begin{equation}
\bigl|(f-[n/n]_f)(z)\bigr|^{1/n}\overset{\mcap}\longrightarrow
e^{-2g^{}_S(z,\infty)},\quad n\to\infty,\quad z\in D.
\label{2.3}
\end{equation}
\end{theoremSta}

In fact, for each $f\in\myA^\circ(\myo\CC\setminus\Sigma)$ with
$\card\Sigma<\infty$, there is only a finite number of the so-called
``spurious'' zero-pole pairs, or Froissart doublets \cite{Fro69}, which makes
impossible the pointwise convergence of PA in Stahl's domain.

The numerical distributions of zeros and poles of PA for the some functions from
Laguerre class are demonstrated on the four
pictures (see Fig. \ref{Fig_pade10_2500_130},
\ref{Fig_pade103_5000_267_full},
\ref{Fig_Pade_6(3000)300_full},
\ref{Fig_PA_t_log_1_300(2000)_full}).

\subsection{Multipoint Pad\'e approximants: Buslaev's Theory}\label{s2s5}
Let the set $\Sigma$ with $\card \Sigma<\infty$, the points $z_k\in \CC\setminus\Sigma$ and functions $f_k\in\myA^\circ(\myo\CC\setminus\Sigma),\, k=1,\dots,m$ be given. We assume that $f_j\in\HH(z_j)$, $j=1,\dots,m$.
Let $n\in\NN$ be fixed.
Then there exists two polynomials $P_n,Q_n\not\equiv0$ of degrees $\leq{n}$
each and such that the following characteristic relations
 \begin{equation}
(Q_nf_j-P_n)(z_j)=O\bigl((z-z_j)^{n_j}\bigr),
\quad z\to z_j,\quad j=1,\dots,m,
\label{bus1}
\end{equation} hold,
where $\sum_{j=1}^m n_j=2n+1$, $n_j\in\ZZ_{+}$, $j=1,\dots,m$. Such
polynomials $P_n$ and $Q_n$ are not unique, but the rational function
$B_n=P_n/Q_n$ is uniquely determined by the relation \eqref{bus1} and
is called a multipoint (or $m$-point) PA of the given set
$\myf=\{f_1,\dots,f_m\}$ of the analytic functions
$f_j\in\myA^\circ(\myo\CC\setminus\Sigma)$. In short, we will call the set
$\myf=\{(f_1,z_1),\dots,(f_m,z_m)\}$ of $m$ multivalued analytic functions
$f_j\in\HH(z_j)$ the {\it multi-germ} or $m$-{\it germ} $\myf$.

In general, all functions of the $m$-germ $\myf$ are supposed to be
different, i.e. not even one of them, say $f_j$, might be obtained as an
analytic continuation of another germ, say $f_k$, $k\neq j$, along paths,
avoiding the set $\Sigma$.

In the generic case,  \eqref{bus1} is equivalent to the relation
\begin{equation}
(f_j-B_n)(z)=O\bigl((z-z_j)^{n_j}\bigr),
\quad z\to z_j,\quad j=1,\dots,m.
\label{bus2}
\end{equation}

We now suppose that in \eqref{bus1} $n_j/n\to 2p_j$ as $n\to\infty$,
$\sum_{j=1}^mp_j=1$, $p_j\geq0$, $j=1,\dots,m$. According to Buslaev's Theory
(2013--2015; see \cite{Bus13}--\cite{Bus15} and
also \cite{Bus15b}, \cite{Bus15c}), there exists (in the nondegenerate case) a
unique (up to a set of zero capacity) compact set $F=F_{\bus}$
which is an $S$-curve weighted in the presence of the external field,
which is generated by the unit negative charge $-\nu$, $\nu=\sum_{j=1}^m
p_j\delta_{z_j}$ concentrated at the interpolation points
$z_1,\dots,z_m$. This compact set possesses the following properties:
$F$ consists of a finite number of analytic arcs, the complement
$\myo\CC\setminus F$ of $F$ consists of a finite number of domains
$D_j\ni z_j$, $\myo\CC\setminus F=\bigcup_{j=1}^mD_j$; each of the
functions $f_j$ is holomorphic (i.e. analytic and single-valued)
in the corresponding domain $D_j$, $f_j\in\HH(D_j)$; if for some $k\neq
j$ the domains coincide with each other, $D_j=D_k$, then the corresponding
functions
are also equal, $f_k=f_j$; the compact set $F$ possesses the property of
``symmetry'' in the external field $V_{*}^{-\nu}$. Namely, the following
relation holds
\begin{equation}
\frac{\partial(V_{}^{\beta_F}-V_{*}^\nu)}{\partial n^{+}}(z)
=
\frac{\partial(V_{}^{\beta_F}-V_{*}^\nu)}{\partial n^{-}}(z),
\quad z\in{F}^\circ,
\label{bus3}
\end{equation}
where $\beta_{F}\in M_1({F})$ is a unique equilibrium probability
measure concentrated on ${F}$ and weighted in $V_{*}^{-\nu}$. In other words, the identity
$$
V_{}^{\beta_{F}}(z)-V_{*}^{\nu}(z)\equiv\const=w_{F},
\quad z\in{F},
$$ is valid,
where ${F}^\circ$ is the union of all open arcs which closures constitute the
compact set $F$; $\partial/\partial n^{\pm}$ are the normal derivatives
to ${F}$ at the point $z\in{F}^\circ$ from the opposite sides of $F$. It is
worth noting that for the fixed $m$-germ $\myf$ the compact set $F$
depends on the numbers $p_j\geq0$, $\sum_{j=1}^m p_j=1$. Therefore,
 the ``optimal'' (Buslaev's) partition of the Riemann sphere into
domains $D_j$ also depends on $p_j$.

Just as in Stahl's Theory, the existence of the $V^{-\nu}$-weighted $S$-curve
$F=F_{\bus}$ is crucial for Buslaev's Theory. In accordance with the theory of Stahl, the weighted $S$-property of the compact set $F$ \eqref{bus3}
may be expressed in the following way
\begin{equation}
\frac{\partial\biggl(\sum\limits_{j=1}^mp_jg_{D_j}(z,z_j)\biggr)}
{\partial n^{+}}
=
\frac{\partial\biggl(\sum\limits_{j=1}^mp_jg_{D_j}(z,z_j)\biggr)}
{\partial n^{-}},\quad z\in F^\circ,
\label{bus3.2}
\end{equation}
where $g_{D_j}(z,z_j)$ is the Green's function for the domain $D_j$ (as
usual, we set $g_{D_j}(z,z_j)\equiv0$ when $z\in D_k\neq D_j$).

In what follows, for the sake of simplicity, we restrict our attention to the
particular case $m=2$ of Buslaev Theorem. Thus, we will discuss
in details only the case of two-point Pad\'e approximant.

Let $z_1=0$, $z_2=\infty$ and $\myf=\{f_0,f_\infty\}$ be the set of two
multivalued analytic functions, such that $f_0\in\HH(0)$ and
$f_\infty\in\HH(\infty)$, and also
$f_0,f_\infty\in\myA^\circ(\myo\CC\setminus\Sigma)$, where
$\card\Sigma<\infty$. Thus, each of the functions $f_0$ and $f_\infty$ is a
multivalued analytic function on the Riemann sphere, punctured at a
finite set of points, each of which is a branch point of $f_0$ or
of $f_\infty$ or of both of them. In other words, $f_0$ and $f_\infty$ are
two germs of the multivalued analytic function, given at the point
$z_1=0$ and $z_2=\infty$, respectively. It is worth noting that they may be
the two germs of the same analytic function, taken at two different points,
namely at $z_1=0$ and $z_2=\infty$.

The two-point (in the classical terminology, this is the $n$-th
truncated fraction of the classical $T$-fraction) PA is defined as follows.
Given a number $n\in\NN$, let $P_n,Q_n\in\CC_n[z]$, $Q_n\not\equiv0$, be polynomials
of degree $\leq{n}$, such that\footnote{For a fixed $n\in\NN$,
we can also claim that the left side of \eqref{b5} is $O(z^{n+1})$
as $z\to0$ and $O(1)$ as $z\to\infty$, but this does not change the
convergence theorem itself.} the following relations hold
\begin{equation}
R_n(z):=\bigl(Q_nf-P_n\bigr)(z)=\begin{cases}
O(z^n),& z\to0,\\
O(1/z),& z\to\infty.%
\end{cases}
\label{b5}
\end{equation}
The pair of polynomials $P_n$ and $Q_n$ is not unique, but the rational
function $B_n:=P_n/Q_n$ is uniquely determined by \eqref{b5},
and is called the two-point diagonal PA of the set of $2$-germ of the
functions $\myf=\{f_0,f_\infty\}$.
In the generic case, it follows from \eqref{b5} that
\begin{equation}
\bigl(f-B_n\bigr)(z)=\begin{cases} O(z^{n}),&z\to0,\\
O(1/z^{n+1}),& z\to\infty.
\end{cases}
\label{b6}
\end{equation}
If it exists, then the rational function $B_n=B_n(z;f)\in\CC_n(z)$ is uniquely
determined by the relation \eqref{b6}.

As for the classical Stahl's case, the existence of an $S$-curve,
associated with the two-point PA and weighted in the external field $V_*^{-\nu}$,
$\nu=(\delta_0+\delta_\infty)/2$, is the crucial element of
Buslaev's two-point convergence theorem. Such a weighted $S$-curve
${F}={F}_{\bus}(f_0,f_\infty)$
exists\footnote{In general there may exist some degenerated cases.}
and realizes the ``optimal'' partition of the Riemann
sphere into two domains $D_0\ni0$ and $D_\infty\ni\infty$, such that
$\myo\CC=D_0\sqcup{F}\sqcup D_\infty$, $f_0\in\HH(D_0)$ and
$f_\infty\in\HH(D_\infty)$. The compact set ${F}$ is a weighted $S$-curve, i.e.
${F}$ consists of a finite number of analytic arcs and possesses the
following property of ``symmetry'':
\begin{equation}
\frac{\partial(V^\beta-V_{*}^\nu)}{\partial n^{+}}(z)
=\frac{\partial(V^\beta-V_{*}^\nu)}{\partial n^{-}}(z),
\quad z\in {F}^\circ,
\label{b7}
\end{equation}
where $\beta=\beta_{F}$ is the probability measure
concentrated on ${F}$ and the equilibrium measure in the external field
$V_{*}^{-\nu}(z)=\frac12\log|z|$, that is,
\begin{equation}
V^\beta(z)-V_{*}^\nu(z)\equiv\const,\quad z\in{F}
\label{b8}
\end{equation}
(In fact, the equilibrium measure $\beta$ is generated by the negative unit charge
$-\nu$, $\nu=(\delta_0+\delta_\infty)/2$).
As before, ${F}^\circ$ is the union of all open arcs of ${F}$ (the
closures of which constitute ${F}$) and $\partial n^{+}$ and $\partial n^{-}$
are the inner (with respect to $D_0$ and $D_\infty$) normal
derivatives at a point $z\in {F}^\circ$ from the opposite sides of
${F}^\circ$. Clearly, $\beta$ is the balayage of the measure $\nu$ from
$D_0\sqcup D_\infty$ onto ${F}$. It is worth noting that ${F}$
itself is a union of the closures of the critical trajectories of a
quadratic differential and the weighted equilibrium measure $\beta=\beta_F$
is given by (see \cite{BuMaSu12})
\begin{equation}
d\beta(\zeta)=\frac1{2\pi i}
\frac1\zeta\sqrt{\frac{V_p(\zeta)}{A_p(\zeta)}}\,d\zeta>0, \quad\zeta\in F.
\label{bus_eq}
\end{equation}

Here, for the sake of simplicity, we only consider the case of two-point PA,
and we set $z_1=0$ and $z_2=\infty$. In what follows, we also suppose that
$f_0$ and $f_\infty$ are the germs of the same multivalued analytic function
$f$, and we denote them by $f_0\in\HH(0)$ and $f_\infty\in\HH(\infty)$.
We suppose that the function $f$ has a finite set of singular points in
$\myo\CC$.

Notice that the functions $f_0(z)=(1-z^2)^{-1/2}\sim1$, $z\to0$, and
$f_\infty=(z^2-1)^{-1/2}\sim1/z$, $z\to\infty$, are the germs of the
same analytic function $f$, given by the equation $(z^2-1)w^2=1$. But the
functions $f_0(z)=(1-z^2)^{-1/2}$ and $f_\infty=(z^2-1)^{-1/2}+1$ are
not so. Thus, the latter case is the generic case, and hence $D_0\cap
D_\infty=\varnothing$ (see
Fig. \ref{Fig_bus210b_4000_120_full}, \ref{Fig_bus205c(2000)195_full}).

Now we are ready to formulate the particular case of Buslaev
Theorem for two-point PA (cf. Stahl Theorem).

\begin{theoremBusTwoPoint}[V. I. Buslaev, 2013--2015]
Let the function $f\in\HH(0)\cap\HH(\infty)$,
$f\in\myA^\circ(\myo\CC\setminus\Sigma)$, $\card\Sigma<\infty$, and let the pair
of germs $f_0,f_\infty$ be in a general
position\footnote{Equivalently, we say that Buslaev's
$S$-curve ${F}$ divides the Riemann sphere into two domains.}.
Let $D_0\sqcup{F}\sqcup D_\infty=\myo\CC$ be the optimal partition of
the Riemann sphere into two domains $D_0\ni0$ and $D_\infty\ni\infty$,
such that $f_0\in\HH(D_0)$, $f_\infty\in\HH(D_\infty)$,
$D_0\cap D_\infty=\varnothing$, and ${F}$
possesses the weighted $S$-property with respect to the external field
$V_*^{-\nu}$, $\nu=(\delta_0+\delta_\infty)/2$. Then for the $n$-diagonal
two-point PA $B_n$ of the set of the germs $\myf=\{f_0,f_\infty\}$ the following
statements hold true:

1) there exists a limit zero-pole distribution for $B_n$, namely,
\begin{equation}
\frac1n\chi(P_n),\frac1n\chi(Q_n)\overset{*}\longrightarrow\beta_F,
\quad n\to\infty;
\label{b9}
\end{equation}

2) there is a convergence in capacity as $n\to\infty$, namely,
\begin{equation}
B_n(z)\overset{\mcap}\longrightarrow f_0(z),\quad z\in D_0,\quad
B_n(z)\overset{\mcap}\longrightarrow f_\infty(z),\quad z\in D_\infty;
\label{b10}
\end{equation}

3) the rate of the convergence in \eqref{b10} is completely
characterized by the relations
\begin{equation}
\begin{aligned}
\bigl|f_0(z)-B_n(z)\bigr|^{1/n}
&\overset{\mcap}\longrightarrow
e^{-g^{}_{D_0}(z,0)},\quad z\in D_0,\\
\bigl|f_\infty(z)-B_n(z)\bigr|^{1/n}
&\overset{\mcap}\longrightarrow
e^{-g^{}_{D_\infty}(z,\infty)},\quad z\in D_\infty.
\end{aligned}
\label{b11}
\end{equation}
\end{theoremBusTwoPoint}

\section{Hermite--Pad\'e polynomials and Hermite approximants}\label{s3}

\subsection{Definition and uniqueness of Hermite approximants}\label{s3s1}
Let us now suppose that the functions $1,f,f^2$ are rationally independent
and let us consider type I HP polynomials, i.e.
$Q_{n,0},Q_{n,1},Q_{n,2}\in\CC_n^*[z]$ and
\begin{equation}
\label{3}
(Q_{n,0}\cdot1+Q_{n,1} \cdot f+Q_{n,2} \cdot f^2)(z)=O\(\frac1{z^{2n+2}}\),
\quad z\to\infty.
\end{equation}
We are now facing two very natural questions. {\it What kind of new results come out from Hermite--Pad\'e
polynomials}? {\it What can be said about the ratios
$Q_{n,0}/Q_{n,2}$ and $Q_{n,1}/Q_{n,2}$ (cf. \eqref{1.4p}), do they converge to
analytic functions corresponding with the given $f$ in some way, or do they
not}? If yes, then does the sequence $\fH_{n,0}(z):=-Q_{n,0}/Q_{n,2}$ provide
more detailed information about the analytic properties of the function $f$ than the
sequence of Pad\'e approximants $[n/n]_f(z)=-P_{n,0}/P_{n,1}$? In
general, the answer is unknown. However, in some special cases the answer is
positive and appears to be very unusual for the HP polynomials theory. Hence, this
problem is very promising for forthcoming investigations.

\begin{lemma}\label{lem1}
Let two triples of polynomials $Q_{n,0}$, $Q_{n,1}$, $Q_{n,2}\in\CC_n^{*}[z]$
and $\myt{Q}_{n,0}$, $\myt{Q}_{n,1}$, $\myt{Q}_{n,2}\in\CC_n^{*}[z]$ satisfy
 relation \eqref{3}. Then the following equalities
\begin{equation}
\label{3.2}
\frac{Q_{n,0}}{\myt{Q}_{n,0}}(z)\equiv
\frac{Q_{n,1}}{\myt{Q}_{n,1}}(z)\equiv
\frac{Q_{n,2}}{\myt{Q}_{n,2}}(z).
\end{equation}
are true.
\end{lemma}

\begin{proof}[Proof of Lemma \ref{lem1}]
Indeed, the conditions of Lemma \ref{lem1} yield
\begin{equation}
\label{3.1}
(\myt Q_{n,0}\cdot1+\myt Q_{n,1} \cdot f+\myt Q_{n,2} \cdot f^2)(z)=O\(\frac1{z^{2n+2}}\),
\quad z\to\infty.
\end{equation}
After multiplying both sides of \eqref{3} by $\myt{Q}_{n,2}$ and both
sides of \eqref{3.1} by ${Q}_{n,2}$, respectively and subtracting the
new equations, we come to
\begin{align}
\label{3.3}
(Q_{n,0}\myt{Q}_{n,2}-\myt{Q}_{n,0}Q_{n,2})(z)
+(Q_{n,1}\myt{Q}_{n,2}-\myt{Q}_{n,1}Q_{n,2})(z)f(z)= \nonumber\\
=O\(\frac1{z^{n+2}}\), \quad z\to\infty.
\end{align}
Just in the same way we obtain the equality
\begin{align}
\label{3.4}
(\myt{Q}_{n,0}Q_{n,1}-Q_{n,0}\myt{Q}_{n,1})(z)
+(Q_{n,1}\myt{Q}_{n,2}-\myt{Q}_{n,1}Q_{n,2})(z)f^2(z)= \nonumber\\
=O\(\frac1{z^{n+2}}\), \quad z\to\infty.
\end{align}
It follows immediately from \eqref{3.3} and \eqref{3.4} that the
polynomial
\[ P_{2n}:=(Q_{n,1}\myt{Q}_{n,2}-\myt{Q}_{n,1}Q_{n,2})\in\CC_{2n}[z] , \]
being of degree $\leq{2n}$, is in fact a type II HP polynomial for the pair
$f,f^2$. Since under the conditions of Lemma \ref{lem1} the triple
$1,f,f^2$ is rationally independent over the field $\CC(z)$, it follows
that in both relations \eqref{3.3} and \eqref{3.4} the order of
approximation at the infinity point should be $O(1/z^{n+1})$ and not
$O(1/z^{n+2})$, unless $P_{2n}\equiv0$. Lemma \ref{lem1} is proved.
\end{proof}

\begin{definition}\label{def1}
In what follows, we call the uniquely defined rational functions
$Q_{n,0}/Q_{n,2}$ and $Q_{n,1}/Q_{n,2}$ the {\it Hermite Approximants}
(HA) $\fH_{n,0}$ and $\fH_{n,1}$, respectively.
\end{definition}


\subsection{Some theoretical results about Hermite approximants}\label{s3s2}
Suppose that $f\in\sL$. Let $Q_{nj},\, j=1,2,3$ be the HP polynomials for the collection $[1,f,f^2]$, and
$\fH_{n,0},\fH_{n,1}$ be the corresponding HA of the function~$f$.

The case (see  \eqref{l1}) $p=2$ and $a_1=-1$, $a_2=1$,
$$
f(z)=\(\frac{z+1}{z-1}\)^\alpha,\quad f(\infty)=1,
$$
where $2\alpha\in\RR\setminus\ZZ$, was treated by
A. Mart\'{\i}nez-Finkelshtein, E. A. Rakhmanov and S. P. Suetin, 2014--2015
(see \cite{MaRaSu13}, \cite{MaRaSu15}).
It was proven \cite[Theorem 1.8]{MaRaSu15} that for $z\in\CC\setminus{F}$ and
$F:=\myo\RR\setminus[-1,1]$, we have for $n\to\infty$ (cf. \eqref{co1}
and \eqref{co2})
\begin{equation}
\begin{aligned}
\frac{Q_{n,1}}{Q_{n,2}}(z)&\to-2\cos{\alpha\pi}
\(\frac{1+z}{1-z}\)^{\alpha},\quad z\notin{F},
\\
\frac{Q_{n,0}}{Q_{n,2}}(z)&\to\(\frac{1+z}{1-z}\)^{2\alpha}=f^2(z),
\quad z\notin{F},\quad f(0)=1.
\end{aligned}
\label{marasu1}
\end{equation}


Let now $f\in\sL$ be given by the representation
\begin{equation}
\label{4}
f(z)=\prod_{j=1}^q\(\frac{z-e_{2j-1}}{z-e_{2j}}\)^\alpha
=\(\prod_{j=1}^q\frac{z-e_{2j-1}}{z-e_{2j}}\)^\alpha,
\quad f(\infty)=1,\quad\alpha\in\RR\setminus\ZZ,
\end{equation}
with $e_j\in\RR$, $-1=e_1<\dots<e_{2q}=1$. We set $\sL_{\RR}$ for this subclass of
$\sL$. Notice that for $f\in\sL_{\RR}$ the pair $f,f^2$
forms the so-called {\it Nikishin's system} (see \cite{Nik86},
\cite{NiSo88}, \cite{Gon03},
\cite{FiLo11}, \cite{ApKu11}, \cite{Lap15}).

Set
$E:=\bigsqcup_{j=1}^q[e_{2j-1},e_{2j}]$,
$D:=\myo\CC\setminus{E}$. Since $E=S$ is the Stahl's compact set for the function  $f$ under consideration,
then by Stahl's Theorem
\begin{equation}
\label{5}
[n/n]_f(z)
\overset{\mcap}\longrightarrow f(z),\quad n\to\infty,
\quad z\in D,
\end{equation}
and
\begin{equation}
\label{6}
\bigl|f(z)-[n/n]_f(z)\bigr|^{1/n}
\overset{\mcap}\longrightarrow e^{-2g_E(z,\infty)}
=e^{2\(\gamma_E-V^\lambda(z)\)},\quad n\to\infty,\quad
z\in D,
\end{equation}
where $g_E(z,\infty)\equiv\gamma_E-V^\lambda(z)$ is the Green's function of $D$,
$\lambda=\lambda_E$ is the unique equilibrium measure of $E$, i.e.
$V^\lambda(x)\equiv\const$, $x\in E$.

Let now $f_2(z)=\const\cdot f(z)$, $\const\neq0$, be another
``branch'' (see \cite{ChCh84}) of the function $f$, which is holomorphic in the domain
$G:=\myo\CC\setminus F$, where $F:=\myo\RR\setminus E$, that is, $G\neq
D$. In general, if $f\in\sL$ is given by the equality
$f(z)=\prod\limits_{j=1}^p(z-a_j)^{\alpha_j}$, then
both functions $f_1=f$ and $f_2$ solve the same differential equation
$$
A_{p}(z)w'+B_{p-2}(z)w=0,
$$
where
\[ A_{p}(z)=\prod_{j=1}^{p}(z-a_j) \ \text{ and } \ B_{p-2}(z)=-A_p(z)\sum_{j=1}^p\alpha_j(z-a_j)^{-1} \]
are polynomials of
degrees $p$ and $p-2$, respectively. If $p=2$, $a_1=-1$, $a_2=1$,
$$
f(z):=\(\frac{z+1}{z-1}\)^\alpha,\quad z\notin E=[-1,1],\quad f(\infty)=1,
$$
then we have
$$
f_2(z)=-2\cos\alpha\pi\(\frac{1+z}{1-z}\)^\alpha,
\quad z\notin F,\quad f_2(0)=1;
$$
see A. Mart\'{\i}nez-Finkelshtein, E. Rakhmanov and S. Suetin \cite{MaRaSu15}.

In wider sense, the following result is valid \cite{Sue15} (cf. \cite{LoMeFi15}).

\begin{theorem}[(S. Suetin, 2015)]\label{the1}
Let $f$ be of type \eqref{4} where $\alpha\in(-1/2,1/2)$, $\alpha\neq0$,
$-1=e_1<\dots<e_{2{q}}=1$. Then

1) all zeros of $Q_{n,0}$, $Q_{n,1}$ and $Q_{n,2}$, up to a finite number
that is fixed and independent of $n$, belong to $F$; there exists a LZD of HP
$Q_{n,j}$:
\begin{gather}
\frac1n\chi(Q_{n,j})\overset{*}\longrightarrow\eta^{\mypha}_F,\quad n\to\infty,
\label{lim_Q}\\
\text{where}
\quad 3V_*^{\eta^{\mypha}_F}(y)+G^{\eta^{\mypha}_F}_E(y)+3g_E(y,\infty)
\equiv\const,\quad y\in F,\quad\supp\eta^{\mypha}_F=F;
\label{pro_1}
\end{gather}

2) the rational function $\fH_{n,1}:=-Q_{n,1}/Q_{n,2}$ interpolates the
function $f_2$ at least at $2n-m$ distinct (``free'') nodes $x_{n,j}$ of
$E^\circ:=\bigsqcup_{j=1}^{q}(e_{2j-1},e_{2j})$ where $m\in\NN$ does not depend
on $n$, and there exist LZD of those free nodes $x_{n,j}$, namely
\begin{gather}
\label{lim_O}
\frac1{2n}\sum_{j=1}^{2n-m}\delta_{x_{n,j}}\overset{*}\longrightarrow\eta^{\mypha}_E,\quad n\to\infty,
\\
\text{where}\quad
3V^{\eta^{\mypha}_E}(x)+G^{\eta^{\mypha}_E}_F(x)\equiv\const,\quad x\in E,
\quad\supp\eta^{\mypha}_E=E;
\label{pro_2}
\end{gather}

3) in the domain $G:=\myo\CC\setminus{F}$, the following relation is valid
\begin{equation}
\fH_{n,1}(z)\overset{\mcap}\longrightarrow
f_2(z),\quad z\in G,\quad n\to\infty;
\label{4.3}
\end{equation}
and the rate of convergence is completely characterized by the relations
(cf. \eqref{6})
\begin{align}
\bigl|f_2(z)-\fH_{n,1}(z)\bigr|^{1/n}
&\overset{\mcap}\longrightarrow
e^{-2G_F^{\eta^{\mypha}_E}(z)}<1,\quad z\in G\setminus{\RR},\quad n\to\infty,
\label{4.4}\\
\varlimsup_{n\to\infty}\bigl|f_2(x)-\fH_{n,1}(x)\bigr|^{1/n}
&\leq e^{-2G_F^{\eta^{\mypha}_E}(x)}<1,\quad x\in E^\circ,
\label{4.4.2}
\end{align}
where the measure $\eta^{\mypha}_E$ solves problem \eqref{pro_2}.
\end{theorem}

In Theorem \ref{the1}
$$
G_F^{\eta^{\mypha}_E}(z)=\int_E g_E(x,z)\,d\eta_E(x)
$$
is the Green potential of the measure $\eta_E$, $\supp\eta_E=E$,
$g_E(x,z)$ is the Green function for $D:=\myo\CC\setminus{E}$,
$$
G_E^{\eta^{\mypha}_F}(z)=\int_F g_F(x,z)\,d\eta_F(x)
$$
is Green potential of measure $\eta_F$, $\supp\eta_F\subset F$,
$g_F(x,z)$ is the Green function for $G:=\myo\CC\setminus{F}$.

Notice that the equilibrium problem \eqref{pro_1} was introduced by S. Suetin
and E. Rakhmanov in \cite{RaSu13} (see also \cite{Sue13b}, \cite{BuSu14b}, \cite{BuSu15})
and is different from the problem that was studied before in
papers \cite{GoRa81}, \cite{Nik86}, \cite{GoRa87}, \cite{GoRaSo97};
see also \cite{NiSo88} and \cite{Gon03}.

The case when we have \eqref{4} with $q=1$ and $e_1=-1$, $e_2=1$, that is,
$$
f(z)=\(\frac{z+1}{z-1}\)^\alpha,\quad f(\infty)=1,
$$
where $2\alpha\in\CC\setminus\ZZ$, was treated by
A. Mart\'{\i}nez-Finkelshtein, E. A. Rakhmanov and S. P. Suetin in 2013--2015.
The first version of Theorem \ref{the1} was established in  \cite[Theorem 1.8]{MaRaSu15}; furthermore, the following explicit
representation for both measures $\eta_E$ and $\eta_F$ were found, namely
\begin{align}
\frac{d\eta_F}{dx}(x)&=
\frac{\sqrt{3}}{2\pi}\frac1{\sqrt[3]{x^2-1}}
\biggl(\frac1{\sqrt[3]{|x|-1}}-\frac1{\sqrt[3]{|x|+1}}\biggr),
\quad x\in\myo\RR\setminus[-1,1],
\notag\\
\frac{d\eta_E}{dx}(x)&
=\frac{\sqrt{3}}{4\pi}\frac1{\sqrt[3]{1-x^2}}
\biggl(\frac1{\sqrt[3]{1-x}}+\frac1{\sqrt{1+x}}\biggr),\quad x\in(-1,1).
\notag
\end{align}
Recall the explicit representation of Chebysh\"ev--Robin equilibrium probability
measure $\lambda_{\Cheb}$ for the unit segment $[-1,1]$:
$$
\frac{d\lambda_{\Cheb}}{dx}=\frac1{\pi\sqrt{1-x^2}},\quad x\in(-1,1).
$$

Under the condition $\alpha=1/3$, i.e. for the function
$$
f(z)=\(\frac{z+1}{z-1}\)^{1/3}
$$
relation \eqref{4.4.2} from Theorem \ref{the1} might be improved in
the following form. The Hermite approximation $\fH_{n,1}(z):=-Q_{n,1}/Q_{n,2}(z)$
possesses the property of ``almost
Chebysh\"ev alternation'' on the open interval $(-1,1)$ in the following
sense. For each positive and arbitrary small $\theta>0$ on the
interval $(-1,1)$ there exist at least $N_n=[2n(1-\theta)]$
consecutive points $x_j$, $-1<x_1<\dots<x_{N_n}<1$, such that the
following equality holds:
\begin{equation}
\label{alt1}
f_2(x_j)-\fH_{n,1}(x_j)
=(-1)^j e^{-2nG_F^{\eta_E}(x_j)}
\biggl\{
\frac23\sqrt[3]{\frac{1+x_j}{1-x_j}\,}\,
\bigl(1+\eps_n(x_j)\bigr)\biggr\},
\end{equation}
where $\eps_n(x)\to0$ as $n\to\infty$ with a geometrical rate locally
uniformly in $(-1,1)$. Let
$$
w_n(z):=e^{2nG_F^{\eta_E}(x_j)}
\frac32\sqrt[3]{\frac{1-x_j}{1+x_j}\,}
$$
be the weight function. Then  \eqref{alt1} implies the following weighted
equality
$$
w_n(x_j)\bigl(f_2(x_j)-\fH_{n,1}(x_j)\bigr)
=(-1)^j (1+\eps_n(x_j)),\quad j=1,\dots,N_n.
$$


\subsection{Orthogonality relations}\label{s3s3}
Let $f\in\HH(\infty)$,
\begin{equation}
\label{n1}
f(z)=\prod_{j=1}^p(z-a_j)^{\alpha_j},
\quad \alpha_j\in\CC\setminus \ZZ,
\quad\sum_{j=1}^p\alpha_j=0,
\end{equation}
where the points $a_j\in\CC$ are pairwise distinct, i.e. $a_j\neq a_k$ when $j\neq
k$. Thus $f\in\myA^\circ(\myo\CC\setminus\Sigma)$, where $\Sigma=\{a_1,\dots,a_p\}$.
 We have in the partial case $f\in\sL_{\RR}$
\begin{equation}
\label{n2}
f(z)=\prod_{j=1}^q\(\frac{z-e_{2j-1}}{z-e_{2j}}\)^{\alpha},
\quad \alpha\in\RR\setminus\ZZ,
\end{equation}
where $-1=e_1<\dots<e_{2q}=1$. Let $|\alpha|\in(0,1/2).$ Let
$E:=\bigsqcup_{j=1}^q[e_{2j-1},e_{2j}]$,
$E^\circ:=\bigsqcup_{j=1}^q(e_{2j-1},e_{2j})$, $E_j:=[e_{2j-1},e_{2j}]$.

We fix the branch of $f$ at $z=\infty$ by $f(\infty)=1$ and
 fix a number $n\in\NN$. By definition \eqref{1.4p}
\begin{equation}
\label{ort1}
\int_\gamma(P_{n,0}+P_{n,1}f)(\zeta)q(\zeta)\,d\zeta=0
\quad\forall q\in\CC_{n-1}[\zeta],
\end{equation}
where $\gamma$ is an arbitrary contour separating the points
$e_1,\dots,e_{2q}$ from the infinity point. Let $f$ be given by \eqref{n2}; then it
follows from \eqref{ort1} that
\begin{equation}
\label{n4}
\int_E P_{n,1}(x)x^k\Delta{f}(x)\,dx=0,\quad k=0,\dots,n-1,
\end{equation}
where $\Delta{f}(x):=f^{+}(x)-f^{-}(x)$, $x\in{E}$.
Since $\const \cdot \Delta{f}>0$ on $E^\circ$ for some $\const\neq0$, we conclude from \eqref{n4}
that:

1) all but some fixed and independent of $n$ number of zeros of $P_{n,1}$
belong to $E$;

2) by Stahl's Theorem, there exists LZD of Pad\'e polynomials $P_{n,1}$:
\begin{equation}
\label{n5}
\frac1n\chi(P_{n,1})\overset{*}\longrightarrow\lambda,
\quad n\to\infty,
\end{equation}
where $\lambda=\lambda_E$ is a unique equilibrium probability measure
concentrated on $E$, i.e.
\begin{equation}
\label{n6}
V^\lambda(x)\equiv\const,\quad x\in E;
\end{equation}
$E=S$ is the Stahl's compact set of $f$.
From definition \eqref{0.2} of HP polynomials, we may write
\begin{equation}
\label{n8}
\int_\gamma (Q_{n,0}+Q_{n,1}f+Q_{n,2}f^2)(\zeta)q(\zeta)\,d\zeta=0
\quad\forall q\in\CC_{2n}[\zeta],
\end{equation}
where $\gamma$ is an arbitrary closed contour that separates points
$e_1,\dots,e_{2q}$ from the infinity point. From \eqref{n8} it follows
that for $q(z)=P_{n+k,1}(z)=P_{n+k,1}(z;f)$ we have
\begin{equation}
\label{n9}
\int_E Q_{n,2}(x)P_{n+k,1}(x)\myt{f}(x)\Delta{f}(x)\,dx=0,
\quad k=1,\dots,n,
\end{equation}
where $\myt{f}(x):=f^{+}(x)+f^{-}(x)$, $x\in{E}$, and
$\const\myt{f}(x)\Delta{f}(x)>0$ for $x\in E^\circ$ with some $\const\neq0$.

From \eqref{n9}, it follows (see \cite{Sue15}) that:

1) all but some fixed and independent of $n$ number of zeros of $Q_{n,2}$
belong to $F:=\myo\RR\setminus{E}$;

2) there exists LZD of HP polynomials $Q_{n,2}$:
$$
\frac1n\chi(Q_{n,2})\overset{*}\longrightarrow\eta_F,
\quad n\to\infty,
$$
where $\eta_F$ is a unique special
equilibrium probability measure concentrated on $F$, i.e.
\begin{equation}
\label{n10}
3V_{*}^{\eta_F}(x)+G_E^{\eta_F}(x)+\psi(x)\equiv\const,
\quad x\in F;
\end{equation}
here
\begin{equation}
\label{n11}
G_E^\mu(z):=\int g_E(\zeta,z)\,d\mu(\zeta),
\quad \psi(z):=3g_E(z,\infty),
\end{equation}
$g_E(\zeta,z)$ is the Green function for $D:=\myo\CC\setminus{E}$. The
pair of compact sets $E,F$ forms the so-called Nuttall condenser
$\sN:=(E,F)=(E;F,\psi)$. We call the corresponding special
equilibrium measure $\eta_F$ from \eqref{n11} the {\it Nuttall
equilibrium measure} (see \cite{RaSu13}, \cite{Sue13b}, \cite{KoSu14}).
For LZD of HP polynomials, the notion of Nuttall's
condenser plays a role, which is very similar to the role played by Stahl's compact set $S$
in the case of Pad\'e polynomials. In general, if the plates
$E,F\not\subset\RR$, then they both possess some special ``symmetry'' property,
see \cite{RaSu13}, \cite{Sue13b}, \cite{KoSu14}.


\subsection{Discussion of some numerical results}\label{s3s4}
We are going to
discuss some numerical examples in order to demonstrate a numerical basis for
Conjectures \ref{conj1} and \ref{conj2} and for
the results of Theorem \ref{the1} as well.

From numerical experiments made by R.~Kovacheva, N.~Ikonomov, and S.~Suetin
\cite{KoIkSu15}, \cite{KoIkSu15b},
it follows that the distribution of zeros of HP polynomials and the convergence
of Hermite approximants itself are very sensitive to the type of branching of
multivalued analytic function. More precisely, the situation
becomes generally much more complicated, even if all
branch points $e_j$ still belong to the real line, but in \eqref{4}
instead of one parameter $\alpha$ we take different parameters
$\alpha_j$, $\alpha_j\in\RR\setminus\ZZ$ (see \eqref{l2}). To be more precise, let the
multivalued analytic function $f$ be given by the explicit
representation
\begin{equation}
f(z)=\prod_{j=1}^{q}\(\frac{z-e_{2j-1}}{z-e_{2j}}\)^{\alpha_j},
\label{2.1-1}
\end{equation}
where $e_1<\dots <e_{2q}$, but $\alpha_j\neq\alpha_k,\,j\neq k$. Let us
fix the germ of $f$ by the relation $f(\infty)=1$.

\begin{mcase}\label{cas1}
Let $q=3$ and
\begin{equation}
\label{f1}
f(z)=
\(\dfrac{z+2.5}{z+1.3}\)^{1/3}\(\dfrac{z+0.8}{z-0.8}\)^{1/3}
\(\dfrac{z-1.3}{z-2.5}\)^{1/3}.
\end{equation}
Since in \eqref{f1} all the exponents are equal to the same
$\alpha=1/3$, the zeros of the associated HP polynomials
$Q_{200,0}$, $Q_{200,1}$, $Q_{200,2}$ of the collection $[1,f,f^2]$
should be distributed in accordance to Theorem \ref{the1}. From
figures \ref{Fig_La1(2000)200_blu_red}--\ref{Fig_La1(2000)200_full}, it
follows that it is really the case. All zeros, except a pair of
Froissart triplets, are distributed on the real line $\RR$ on the complement
of three real segments $[-2.5,-1.3]$, $[-0.8,0.8]$, and $[1.3,2.5]$.
\end{mcase}

\begin{figure}[!ht]
\centerline{
\includegraphics[width=15cm,height=15cm]{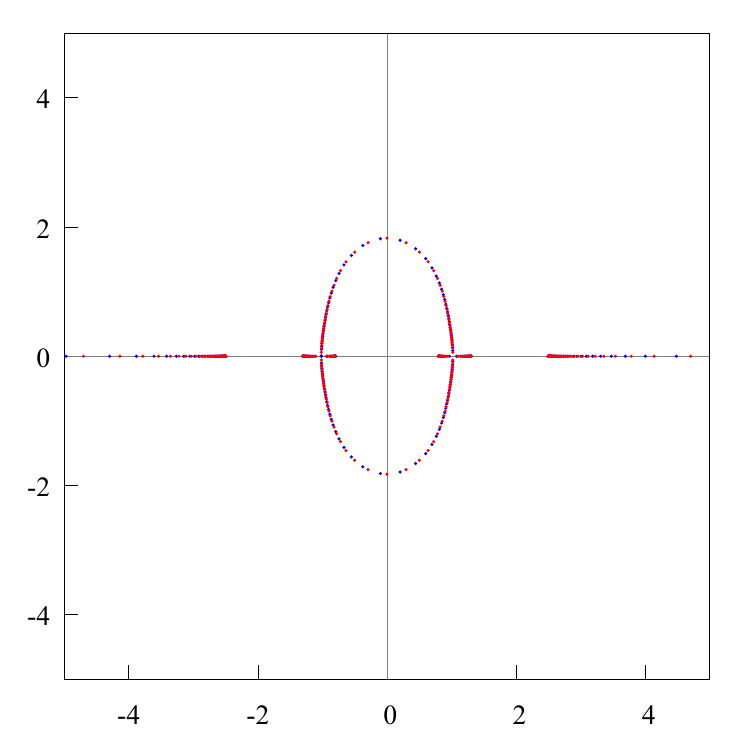}}
\vskip-6mm
\caption{Zeros of HP polynomials $Q_{320,0}$ (blue points) and $Q_{320,1}$ (red
points) for the triple of functions
for $[1,f,f^2]$, where
$f(z)=
\(\dfrac{z+2.5}{z+1.3}\)^{1/3}{\(\dfrac{z+.8}{z-.8}\)^{-1/3}}
\(\dfrac{z-1.3}{z-2.5}\)^{1/3}$.
There is a membrane which splits the complement of the segments
$[-2.5,-1.3]$, $[-0.8,0.8]$ and $[1.3,2.5]$ into two domains.
Both domains are simply connected.
The zeros of these HP
polynomials are distributed on real line $\RR$ on the complement of the segments
$[-2.5,-1.3]$, $[-0.8,0.8]$ and $[1.3,2.5]$ and on
this membrane. The points of intersection of the membrane with the two
segments are the Chebotar\"ev's points of zero-density for the
equilibrium measure for a compact set~$F$.
By chance, there are no Froissart doublets at all;
cf. Fig.~\ref{Fig_La2(1500)320_full}.
}
\label{Fig_La2(1500)320_blu_red}
\end{figure}

\begin{figure}[!ht]
\centerline{
\includegraphics[width=15cm,height=15cm]{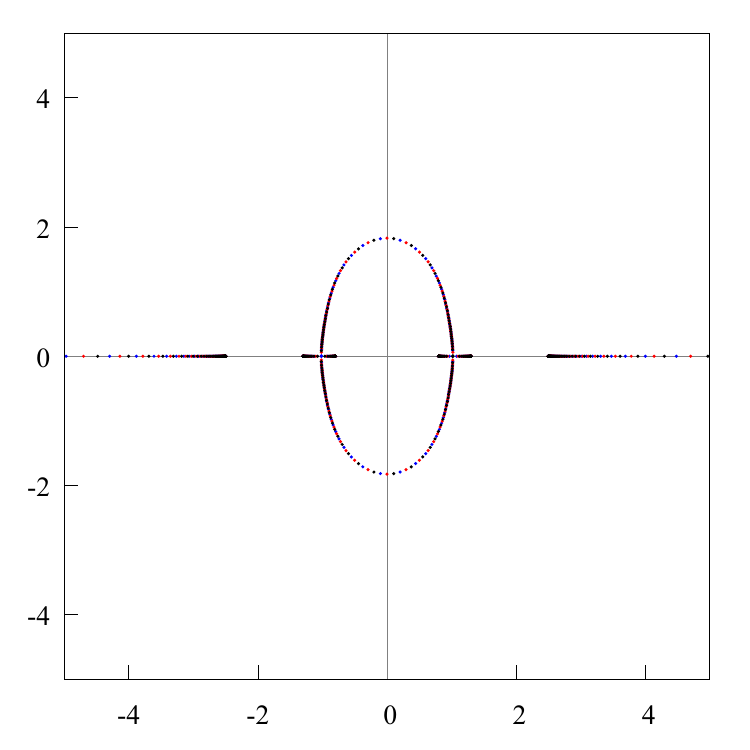}}
\vskip-6mm
\caption{Zeros of HP polynomials $Q_{320,0}$ (blue points), $Q_{320,1}$ (red
points), and $Q_{320,2}$ (black points) for the triple of functions
for $[1,f,f^2]$, where
$f(z)=
\(\dfrac{z+2.5}{z+1.3}\)^{1/3}{\(\dfrac{z+0.8}{z-.08}\)^{-1/3}}
\(\dfrac{z-1.3}{z-2.5}\)^{1/3}$.
There is a membrane which splits the complement to the segments
$[-2.5,-1.3]$, $[-0.8,0.8]$ and $[1.3,2.5]$ into two domains.
Both domains are simply connected.
The zeros of these HP
polynomials are distributed on real line $\RR$ on the complement of the segments
$[-2.5,-1.3]$, $[-0.8,0.8]$ and $[1.3,2.5]$ and on
this membrane. The points of intersection of the membrane with the two
segments are the Chebotar\"ev's points of zero-density for the
equilibrium measure for a compact set $F$.
By chance, there are no Froissart triplets at all;
cf. Fig.~\ref{Fig_La2(1500)320_blu_red}.
}
\label{Fig_La2(1500)320_full}
\end{figure}

\begin{figure}[!ht]
\centerline{
\includegraphics[width=15cm,height=15cm]{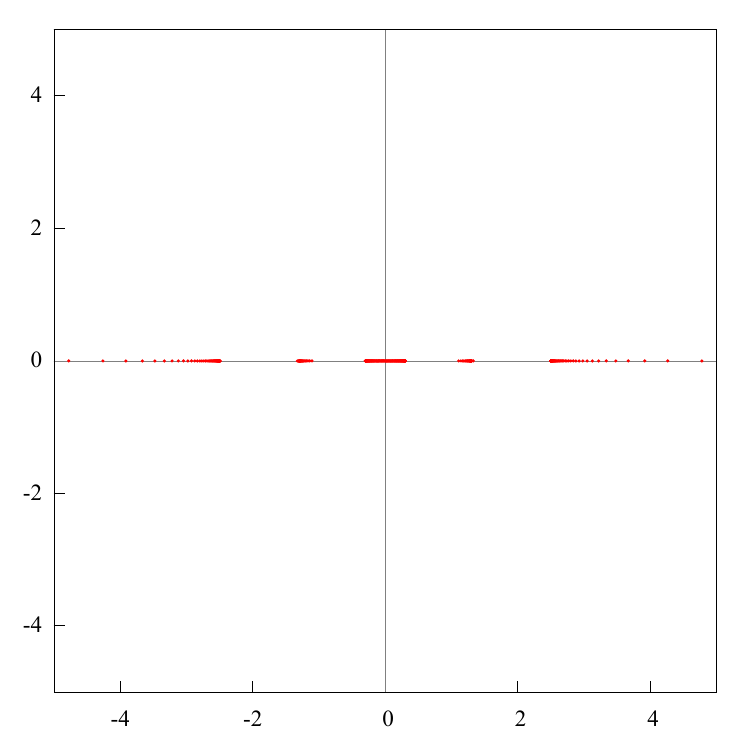}}
\vskip-6mm
\caption{Zeros of HP polynomial $Q_{320,1}$ (red points) for the triple of functions
$[1,f,f^2]$, where
$f(z)=
\(\dfrac{z+2.5}{z+1.3}\)^{1/3}{\(\dfrac{z+0.3}{z-0.3}\)^{1/2}}
\(\dfrac{z-1.3}{z-2.5}\)^{1/3}$.
There is no membrane here.
The zeros of HP polynomial $Q_{320,1}$ are distributed
on the real line on the complement of four segments $[-2.5,-1.3]$, $[-a,-0.3]$,
$[0.3,a]$, and $[1.3,2.5]$ where $a\in(0.3,1.3)$ is an unknown parameter. This
parameter should be evaluated from an appropriate theoretical-potential
equilibrium problem.
By chance, there are no Froissart doublets at all.
In case of the given function $f$, the numerical distribution of zeros of HP
polynomial $Q_{320,1}$ is very
different from numerical distribution of zeros of HP polynomial $Q_{320,0}$
and $Q_{320,2}$;
cf. Fig. \ref{Fig_la1sq(1500)320_blu},
 \ref{Fig_la1sq(1500)320_blk}, \ref{Fig_la1sq(1500)320_full}.
}
\label{Fig_la1sq(1500)320_red}
\end{figure}

\begin{figure}[!ht]
\centerline{
\includegraphics[width=15cm,height=15cm]{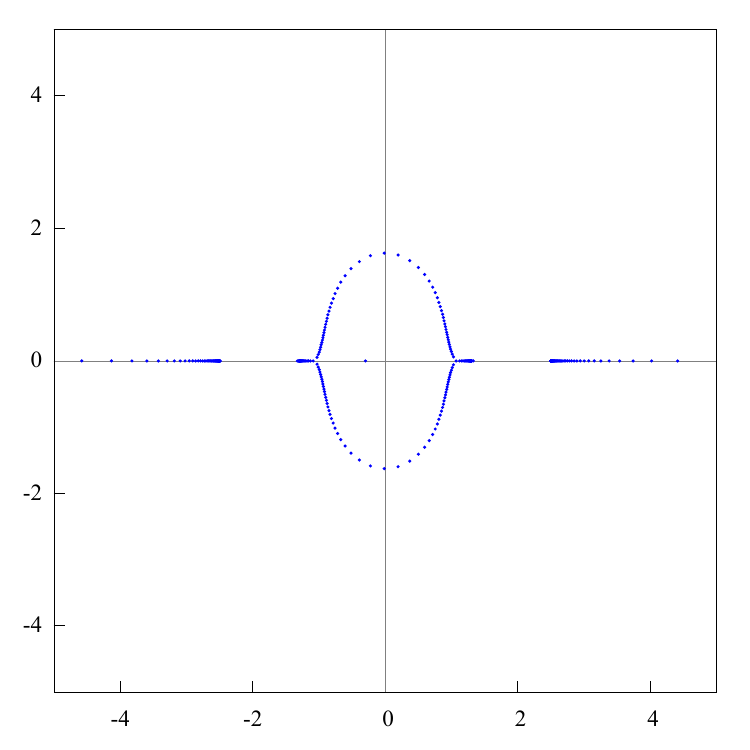}}
\vskip-6mm
\caption{Zeros of HP polynomial $Q_{320,0}$ (blue points)
for the triple of functions $[1,f,f^2]$, where
$f(z)=
\(\dfrac{z+2.5}{z+1.3}\)^{1/3}{\(\dfrac{z+0.3}{z-0.3}\)^{1/2}}
\(\dfrac{z-1.3}{z-2.5}\)^{1/3}$.
There is a membrane which splits the Riemann sphere into two domains.
The zeros of HP polynomial $Q_{320,0}$ are distributed on this membrane
and on the real line on the complement of three segments $[-2.5,-1.3]$,
$[-a,a]$, and $[1.3,2.5]$. The membrane and a parameter $a\in(0.3,1.3)$ come
from an appropriate theoretical-potential equilibrium problem.
By chance, there are no Froissart doublets at all.
In case of the given function $f$, the numerical distribution of zeros of HP
polynomial $Q_{320,0}$ is very
different from the numerical distribution of zeros of HP polynomial $Q_{320,1}$;
cf. Fig. \ref{Fig_la1sq(1500)320_red} and also \ref{Fig_la1sq(1500)320_full}.
There is a zero of polynomial $Q_{320,0}$ inside the membrane which
correspond to the simple zero of the function $f^2$ at the point $z=-0.3$.
}
\label{Fig_la1sq(1500)320_blu}
\end{figure}

\begin{figure}[!ht]
\centerline{
\includegraphics[width=15cm,height=15cm]{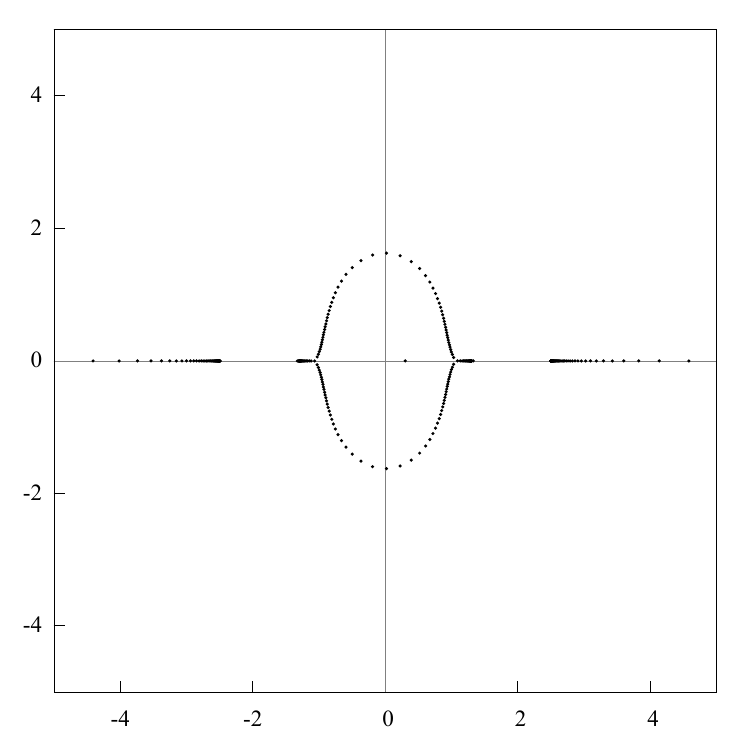}}
\vskip-6mm
\caption{Zeros of HP polynomial $Q_{320,2}$ (black points)
for the triple of functions $[1,f,f^2]$, where
$f(z)=
\(\dfrac{z+2.5}{z+1.3}\)^{1/3}{\(\dfrac{z+0.3}{z-0.3}\)^{1/2}}
\(\dfrac{z-1.3}{z-2.5}\)^{1/3}$.
There is a membrane which splits the Riemann sphere into two domains.
The zeros of HP polynomial $Q_{320,2}$ are distributed on this membrane
and on the real line on the complement of three segments $[-2.5,-1.3]$,
$[-a,a]$, and $[1.3,2.5]$. The membrane and a parameter $a\in(0.3,1.3)$ come
from an appropriate theoretical-potential equilibrium problem.
By chance, there are no Froissart doublets at all.
In case of the given function $f$, the numerical distribution of zeros of HP
polynomial $Q_{320,2}$ is very
different from the numerical distribution of zeros of HP polynomial $Q_{320,1}$;
cf. Fig. \ref{Fig_la1sq(1500)320_red} and also \ref{Fig_la1sq(1500)320_full}.
There is a zero of polynomial $Q_{320,2}$ inside the membrane which
correspond to the simple pole of the function $f^2$ at the point $z=0.3$.
}
\label{Fig_la1sq(1500)320_blk}
\end{figure}

\begin{figure}[!ht]
\centerline{
\includegraphics[width=15cm,height=15cm]{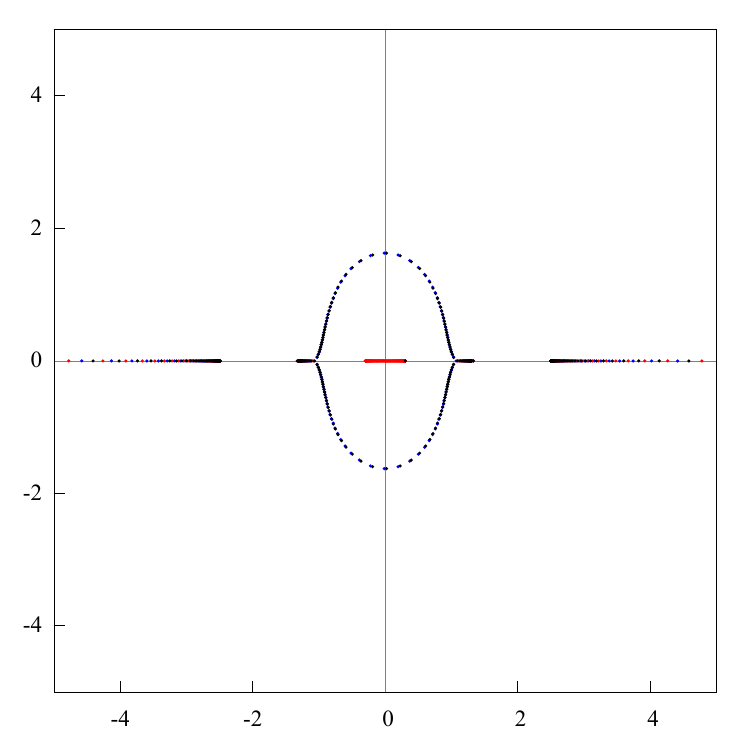}}
\vskip-6mm
\caption{Zeros of HP polynomials $Q_{320,0}$ (blue points),
$Q_{320,1}$ (red points), and $Q_{320,2}$ (black points) for the triple of functions
$[1,f,f^2]$, where
$f(z)=
\(\dfrac{z+2.5}{z+1.3}\)^{1/3}{\(\dfrac{z+0.3}{z-0.3}\)^{1/2}}
\(\dfrac{z-1.3}{z-2.5}\)^{1/3}$.
There is a membrane which splits the Riemann sphere into two domains.
One domain is simply connected, but the other is a doubly connected domain.
The zeros of these HP polynomials are distributed
in accordance with the description given in
Fig. \ref{Fig_la1sq(1500)320_red},
\ref{Fig_la1sq(1500)320_blu}, \ref{Fig_la1sq(1500)320_blk}.
The points of intersection of the membrane with the two
segments $[-1.3,-a]$ and $[a,1.3]$ are the Chebotar\"ev's points of
zero-density for the equilibrium measure
from an appropriate theoretical-potential equilibrium problem.
By chance, there are no Froissart doublets at all.
}
\label{Fig_la1sq(1500)320_full}
\end{figure}

In the general situation \eqref{2.1-1}, when there are different
$\alpha_j$ (instead of a single $\alpha$) there should be membranes which separate
the segments of the set~$F$
(see Fig. \ref{Fig_La2(1500)320_blu_red}--\ref{Fig_la1sq(1500)320_full}).

\begin{mcase}\label{cas2}
Let $q=3$ and
\begin{equation}
\label{f2}
f(z)=
\(\dfrac{z+2.5}{z+1.3}\)^{1/3}\(\dfrac{z+0.8}{z-0.8}\)^{-1/3}
\(\dfrac{z-1.3}{z-2.5}\)^{1/3}.
\end{equation}
Thus in \eqref{2.1-1} $\alpha_j=(-1)^{j+1}\alpha$. Figures
\ref{Fig_La2(1500)320_blu_red}--\ref{Fig_La2(1500)320_full}
represent the numerical distribution of zeros of HP polynomials
$Q_{320,0}$, $Q_{320,1}$, $Q_{320,2}$ of the collection of the functions
$[1,f,f^2]$. In this case there is a
membrane, which splits the complement of the segments $[-2.5,-1.3]$,
$[-0.8,0.8]$ and $[1.3,2.5]$ into two domains.
The zeros of these HP
polynomials are distributed on the real line $\RR$ on the complement of the segments
$[-2.5,-1.3]$, $[-.8,.8]$ and $[1.3,2.5]$ and on
this membrane. The points of intersection of the membrane with the two
segments are the Chebotar\"ev's points of zero-density for the
equilibrium measure for a compact set $F$.
By chance, there are no Froissart triplets at all
(see Fig. \ref{Fig_La2(1500)320_blu_red}--\ref{Fig_La2(1500)320_full}).
\end{mcase}

\begin{mcase}\label{cas3}
Let $q=3$ and
\begin{equation}
\label{f3}
f(z)=
\(\dfrac{z+2.5}{z+1.3}\)^{1/3}\(\dfrac{z+0.3}{z-0.3}\)^{1/2}
\(\dfrac{z-1.3}{z-2.5}\)^{-1/3}.
\end{equation}
Figures \ref{Fig_la1sq(1500)320_red}--\ref{Fig_la1sq(1500)320_full}
represent the numerical distribution of zeros of HP polynomials
$Q_{200,0}$, $Q_{200,1}$, $Q_{200,2}$ for the collections of functions
$[1,f,f^2]$. There also exists a membrane, but of another
type than in Case \ref{cas2}. This membrane splits the complement of
the three segments $[-2.5,-1.3]$, $[-0.3,0.3]$ and $[1.3,2.5]$ into two
domains. The zeros of those HP polynomials are distributed
on the real line $\RR$ on the complement of the segments
$[-2.5,-1.3]$, $[-0.3,0.3]$ and $[1.3,2.5]$ on this new membrane. Just as in Case \ref{cas2},
the two points of intersection of the membrane with the segments are the
Chebotarev's points of zero-density for the equilibrium measure for
compact set~$F$.
\end{mcase}

\subsection{Final remarks about Hermite approximants}\label{s3s5}

Thus, from the numerical experiments of R. Kovacheva, N. Ikonomov, and
S. Suetin, see \cite{KoIkSu15}, \cite{KoIkSu15b}
it follows
that the distribution of zeros of HP polynomials for the collection
$[1,f,f^2]$ and the convergence of Hermite approximants $\fH_{n,j}$, $j=0,1$, itself
are very sensitive to the type of branching of the given multivalued
analytic function $f$. By this reason, it might be very difficult to
construct a general theory of limit zero distribution of HP polynomials
of such type as Stahl's
 and Buslaev's theories are. But as surplus, this sensitivity makes Hermite
approximants $\fH_{n,j}$ very powerful tool to recover the unknown properties
of a multivalued analytic function given by a germ.



\clearpage
\newpage

\clearpage
\newpage

\clearpage
\newpage

\clearpage
\newpage


\clearpage
\newpage

\clearpage
\newpage

\clearpage
\newpage

\clearpage
\newpage

\clearpage
\newpage

\clearpage
\newpage

\clearpage
\newpage

\clearpage
\newpage

\clearpage
\newpage

\clearpage
\newpage

\end{document}